\documentclass[preprint,10pt]{elsarticle}
\usepackage[margin=15mm]{geometry}
\makeatletter
\def\ps@pprintTitle{%
  \let\@oddhead\@empty
  \let\@evenhead\@empty
  \let\@oddfoot\@empty
  \let\@evenfoot\@oddfoot
}
\makeatother

\usepackage{amssymb}
\usepackage{amsmath}
\usepackage{amsthm}
\usepackage{caption}
\usepackage{subcaption}
\usepackage[linesnumbered,ruled,vlined]{algorithm2e}
\usepackage{float}
\usepackage{booktabs}
\usepackage{lineno}
\usepackage{multirow}
\usepackage[table]{xcolor}

\usepackage[colorlinks=true,linkcolor=blue,citecolor=blue,urlcolor=blue]{hyperref}
\usepackage[capitalize,nameinlink]{cleveref}
\usepackage{doi}



\crefname{equation}{equation}{equations}
\Crefname{equation}{Equation}{Equations}
\crefformat{section}{\textbf{#2Section~#1#3}}
\Crefformat{section}{\textbf{#2Section~#1#3}}
\crefformat{subsection}{\textbf{#2Section~#1#3}}
\Crefformat{subsection}{\textbf{#2Section~#1#3}}
\crefformat{theorem}{\textbf{#2Theorem~#1#3}}
\Crefformat{theorem}{\textbf{#2Theorem~#1#3}}
\crefformat{corollary}{\textbf{#2Corollary~#1#3}}
\Crefformat{corollary}{\textbf{#2Corollary~#1#3}}
\crefformat{definition}{\textbf{#2Definition~#1#3}}
\Crefformat{definition}{\textbf{#2Definition~#1#3}}
\crefformat{remark}{\textbf{#2Remark~#1#3}}
\Crefformat{remark}{\textbf{#2Remark~#1#3}}
\crefname{figure}{Figure}{Figures}
\Crefname{figure}{Figure}{Figures}
\crefname{algorithm}{\textbf{Algorithm}}{\textbf{Algorithms}}
\Crefname{algorithm}{\textbf{Algorithm}}{\textbf{Algorithms}}
\crefname{example}{\textbf{Example}}{\textbf{Examples}}
\Crefname{example}{\textbf{Example}}{\textbf{Examples}}

\DeclareCaptionLabelFormat{periodlabel}{#2.}
\usepackage{bm}

\journal{Applied Numerical Mathematics}

\providecommand{\norm}[1]{\left\Vert#1\right\Vert}

\newtheorem{theorem}{Theorem}

\newtheorem{remark}{Remark}

\newtheorem{definition}{Definition}
\newdefinition{example}{Example}
\newproof{pf}{Proof}
\newcommand{\proofofref}{}
\newproof{zproofof}{Proof of \proofofref}

\renewenvironment{example}
  {\begin{oldexample}}
  {\hfill$\blacksquare$\end{oldexample}}

\begin{document}

\begin{frontmatter}



\title{Implementation of Milstein Schemes for Stochastic Delay-Differential Equations with Arbitrary Fixed Delays}


\begin{NoHyper}

\author[Label1]{Mitchell Griggs\corref{cor1}} 
\ead{mitchell.griggs@connect.qut.edu.au}
\author[Label1,Label2]{Kevin Burrage}
\ead{kevin.burrage@qut.edu.au}
\author[Label1]{Pamela Burrage}
\ead{pamela.burrage@qut.edu.au}

\affiliation[Label1]{organization={School of Mathematical Sciences, Queensland University of Technology (QUT)},
            city={Brisbane},
            country={Australia}}
\affiliation[Label2]{organization={Visiting Professor, Department of Computer Science, University of Oxford},
            city={Oxford},
            country={United Kingdom}}
\cortext[cor1]{Corresponding author}

\end{NoHyper}

\begin{abstract}
This paper develops methods for numerically solving stochastic delay-differential equations (SDDEs) with multiple fixed delays that do not align with a uniform time mesh. We focus on numerical schemes of strong convergence orders $1/2$ and $1$, such as the Euler--Maruyama and Milstein schemes, respectively. Although numerical schemes for SDDEs with delays $\tau_1,\ldots,\tau_K$ are theoretically established, their implementations require evaluations at both present times such as $t_n$, and also at delayed times such as $t_n-\tau_k$ and $t_n-\tau_l-\tau_k$. As a result, previous simulations of these schemes have been largely restricted to the case of divisible delays. We develop simulation techniques for the general case of indivisible delays where delayed times such as $t_n-\tau_k$ are not restricted to a uniform time mesh. To achieve order of convergence (OoC) $1/2$, we implement the schemes with a fixed step size while using linear interpolation to approximate delayed scheme values. To achieve OoC $1$, we construct an augmented time mesh that includes all time points required to evaluate the schemes, which necessitates using a varying step size. We also introduce a technique to simulate delayed iterated stochastic integrals on the augmented time mesh, by extending an established method from the divisible-delays setting. We then confirm that the numerical schemes achieve their theoretical convergence orders with computational examples.
\end{abstract}





\end{frontmatter}



\section{Introduction}\label{Section1}

In this paper, we develop practical methods to produce numerical solutions for It\^o SDDEs of the form
\begin{align}
\mathrm{d}X(t)&=[A_0X(t)+f(t,X(t),X(t-\tau_1),\ldots,X(t-\tau_K))]\,\mathrm{d}t\nonumber\\
&\quad+\sum_{j=1}^m[A_jX(t)+g_j(t,X(t),X(t-\tau_1),\ldots,X(t-\tau_K))]\,\mathrm{d}W_j(t),\quad t\in[0,T],\label{Equation1}\\
X(t)&=\phi(t),\quad t\in[-\tau,0].\nonumber
\end{align}
Equation \eqref{Equation1} includes $d$-dimensional solution process $X=(X(t))_{t\in[0,T]}$, constant matrices $A_0,\ldots,A_m\in\mathbb{R}^{d\times d}$, functions $f:[0,T]\times\mathbb{R}^d\times(\mathbb{R}^d)^K\to\mathbb{R}^d$ and $g_1,\ldots,g_m:[0,T]\times\mathbb{R}^d\times(\mathbb{R}^d)^K\to\mathbb{R}^d$, constant delays $\tau_1,\ldots,\tau_K>0$, history process $\phi$, and independent Wiener processes $W_1,\ldots,W_m$. The largest delay is $\tau=\max\{\tau_1,\ldots,\tau_K\}$, the history process satisfies the moment condition $\sup_{t\in[-\tau,0]}\mathbb{E}[\norm{\phi(t)}^4]<\infty$, and each of $f,g_1,\ldots,g_m$ is $C^1$ in the first component ($t$) while also being $C^2$ in all other components ($X(t),X(t-\tau_1),\ldots,X(t-\tau_K)$). Without loss of generality, we suppose that one of the delays divides $T$. That is, for some $k$, there exists $N_k\in\mathbb{N}$ such that $N_k\tau_k=T$. If this is not the case then we simply include another delay, $\tau_{K+1}=T$, and include $X(t-\tau_{K+1})$ in the arguments for $f,g_1,\ldots,g_m$.

Equation \eqref{Equation1} arises in models involving both randomness and time delays. These models have recently found applications in economics \cite{LeeKimKim2011,HodgsonReisinger2022}, medicine \cite{BouletBalasubramaniamDaffertshoferLongtin2010,HartungTuri2013}, ecology \cite{WangWangChen2019}, and other fields. 
 In these and other models, properties of the solution trajectories depend on the delay values. \textcolor{black}{For example, financial options may be priced depending on delayed stochastic volatility.} 

\begin{example}\label{FinancialExample}
\textcolor{black}{The delayed local-volatility model (DLVM) from mathematical finance involves a portfolio of share prices, $X(t)\in\mathbb{R}^d$ at time $t$, modelled using SDDEs whose volatility terms ($g_j$) include delay dependence. Option values typically take the form $V=\mathbb{E}[\nu(X(T))]$ for a payoff function $\nu$. Simpler DLVMs include discrete delays such as in equation \eqref{Equation1}, and also admit analytic expressions for $V$ for particular payoff functions $\nu$ \cite{McWilliams2011}. However, more general DLVMs include distributed delays \cite{GuineaJuliaCaroCarretero2024} such as
\begin{equation}D_{-\tau}^0(X(t))=\int_{-\tau}^0 K(s+\tau)X(t+s)\,\mathrm{d}s,\label{DistributedDelay}\end{equation}
for a kernel $K$. Techniques such as quadrature approximations of \eqref{DistributedDelay} sample $X(t+s)$ at various times $s\in[-\tau,0]$, yielding discretisations with delays that are generally indivisible.}
\end{example}

In our work, we allow for general delays, that may or may not align with a uniform time mesh. We develop simulation techniques to solve \eqref{Equation1} numerically, with convergence orders $1/2$ and $1$. The standard method for numerically solving (non-delayed) stochastic differential equations (SDEs) is the Euler--Maruyama (EM) approximation by Maruyama \cite{Maruyama1955}, which has OoC $1/2$. The EM scheme is easily extendable to \eqref{Equation1}, and only includes additional calculations for delayed scheme values. Achieving OoC $1$ with SDE schemes requires alternative techniques such as the method by Milstein \cite{Milstein1975}, which includes the iterated stochastic integrals
\begin{equation}I_{ij}(t_n,t_{n+1})=\int_{t_n}^{t_{n+1}}\int_{t_n}^s\,\mathrm{d}W_i(u)\,\mathrm{d}W_j(s),\quad i,j=1,\ldots,m.\label{Iij}\end{equation}
The Milstein scheme is further complicated when delays are included (except for the case of delay-free noise, where $g_j(t,X(t),X(t-\tau_1),\ldots,X(t-\tau_K))=g_j(t,X(t))$, $j=1,\ldots,m$), as this includes the delayed iterated stochastic integrals
\begin{equation}I_{ij}^{\tau_k}(t_n,t_{n+1})=\int_{t_n}^{t_{n+1}}\int_{t_n}^s\,\mathrm{d}W_i(u-\tau_k)\,\mathrm{d}W_j(s)=\int_{t_n}^{t_{n+1}}\int_{t_n-\tau_k}^{s-\tau_k}\,\mathrm{d}W_i(u)\,\mathrm{d}W_j(s),\quad i,j=1,\ldots,m,\quad k=1,\ldots,K.\label{Iijtauk}\end{equation}
The SDDE Milstein scheme is given by K\"uchler and Platen \cite{KuchlerPlaten2000} when $K=1$, but simulation methods for \eqref{Iij} were still in development at the time of this research, so examples from this time only include delay-free noise.

Hofmann and M\"uller-Gronbach \cite{HofmannMullerGronbach2006} give an SDDE Milstein scheme for $K$ delays, using the approximations
\begin{equation}I_{ij}(t_n,t_{n+1})\approx\frac{(W_i(t_{n+1})-W_i(t_n))\,(W_j(t_{n+1})-W_j(t_n))}{2}\label{IijSimple}\end{equation}
and
\begin{equation}I_{ij}^{\tau_k}(t_n,t_{n+1})\approx\frac{(W_i(t_{n+1}-\tau_k)-W_i(t_n-\tau_k))\,(W_j(t_{n+1})-W_j(t_n))}{2}\label{IijtaukSimple}\end{equation}
to simulate the iterated integrals. This produces an efficient Milstein approximation of \eqref{Equation1}, but with OoC $1/2$. One of our objectives in this paper is to improve the approximation \eqref{IijtaukSimple} of \eqref{Iijtauk}, to improve the SDDE Milstein scheme so as to successfully achieve demonstrable OoC $1$. To produce an improved approximation for \eqref{Iijtauk}, we find inspiration from techniques that have been developed for the non-delayed case, which we extend to the setting with delays.

One method to approximate \eqref{Iij} is the trapezoidal method by Milstein \cite{Milstein1995}. Cao, Zhang, and Karniadakis \cite{CaoZhangKarniadakis2015} extend the trapezoidal method to simulate \eqref{Iijtauk}, for the SDDE Milstein scheme in the case of a single delay that divides $T$. We extend this approximation further, to the case of multiple delays, including indivisible ones where some delay $\tau_k$ does not divide $T$. Our focus is on implementation techniques, to complement existing theory. For example, there exists a continuous Milstein scheme for \eqref{Equation1}, by Kloeden and Shardlow \cite{KloedenShardlow2012}, but this scheme requires simulation of $I_{ij}^{\tau_k}(t_n,t)$ for every time $t$. However, it is not possible to simulate such a scheme on a continuum of values, due to the fluctuating iterated stochastic integrals. Therefore we restrict simulations to discrete mesh points, simulating the numerical schemes at every point needed for the Milstein scheme. This requires a varying step size when the delays are indivisible. We give the construction of this time mesh specifically for the Milstein scheme, but our implementation method is readily applicable for other schemes, which we demonstrate with the SDDE Magnus--Milstein (MM) method, by Griggs, Burrage, and Burrage \cite{GriggsBurrageBurrage2025}.



\textcolor{black}{Future extensions of this work include solving SDDEs numerically with non-constant or distributed delays. For SDDEs with time-dependent and state-dependent delays of the form $\tau_k=\tau_k(t,X(t))$, numerical schemes of OoC $1/2$ are already well suited, typically relying on approximations of delayed values such as $X(t_n-\tau_k(t_n,X(t_n)))$. Higher-order schemes, however, require further investigation. One possible direction is the development of hybrid approaches that combine mesh augmentation with local approximations of delayed arguments, with the goal of improving accuracy beyond OoC $1/2$. 
 Another generalisation of this paper concerns distributed delays, capturing process history through an integral of the form \eqref{DistributedDelay}. In this setting, an SDDE such as
\begin{equation}\mathrm{d}X(t)=f\left(t,X(t),D_{-\tau}^0(X(t))\right)\,\mathrm{d}t+\sum_{j=1}^mg_j\left(t,X(t),D_{-\tau}^0(X(t))\right)\,\mathrm{d}W_j(t)\label{DistributedSDDE}\end{equation}
may be treated numerically by using a quadrature approximation for $D_{-\tau}^0(X(t))$, reducing the equation to one with a finite collection of delays, $0<\tau_1<\cdots<\tau_K=\tau$. That is, \eqref{DistributedSDDE} may be approximated by an SDDE of the form 
\begin{equation*}
\mathrm{d}X(t)=f\left(t,X(t),X(t-\tau_1),\ldots,X(t-\tau_K)\right)\,\mathrm{d}t+\sum_{j=1}^mg_j\left(t,X(t),X(t-\tau_1),\ldots,X(t-\tau_K)\right)\,\mathrm{d}W_j(t),
\end{equation*}
which is a particular case of \eqref{Equation1}.} 



This paper gives instructions for simulating numerical SDDE schemes with strong OoC $1$, and is presented in the following structure. In \cref{Section2}, we present our selected schemes, which are later used while constructing the augmented time mesh. We discuss implementation complications for these schemes, in \cref{Section3}. We simulate these schemes with OoC $1/2$ in \cref{Section4}, by approximating delayed scheme values with linear interpolants. We begin improving the convergence orders of these schemes in \cref{Section5}, by augmenting the simulation time mesh. We then extend the trapezoidal approximation (for delayed iterated stochastic integrals) to this augmented mesh, in \cref{Section6}. We provide example simulations of the schemes on the augmented mesh, in \cref{Section7}. 
 In \cref{Section8}, we conclude this work with discussion of the size of the augmented mesh, the applicability of our technique for higher-order schemes, and also the significance of our method.

\section{Selected numerical schemes}\label{Section2}

In this section, we define our selected numerical schemes used to approximate solutions of \eqref{Equation1}. We begin by defining a solution and stating an associated existence-uniqueness result. The numerical schemes are derived from applying regular SDE techniques across Bellman intervals, which convert the SDDE into a system of SDEs. This also enables us to define analogous properties such as convergence orders. We then describe simulation techniques for the integrals \eqref{Iij} and \eqref{Iijtauk}, which we use in numerical examples demonstrating the convergence orders of the schemes, for particular (divisible) delays.


\begin{definition}
A stochastic process $X=(X(t))_{t\in[-\tau,T]}$ is a \emph{(strong) solution} of equation \eqref{Equation1} if $X$ satisfies $X(t)=\phi(t)$ when $t\leqslant0$, and
\begin{align*}
X(t)&=X(0)+\int_0^t[A_0X(s)+f(s,X(s),X(s-\tau_1),\ldots,X(s-\tau_K))]\,\mathrm{d}s\\
&\quad\quad+\sum_{j=1}^m\int_0^t[A_jX(s)+g_j(s,X(s),X(s-\tau_1),\ldots,X(s-\tau_K))]\,\mathrm{d}W_j(s),\quad t\in[0,T].
\end{align*}
Solutions of equation \eqref{Equation1} are \emph{(pathwise) unique} if whenever $X=(X(t))_{t\in[-\tau,T]}$ and $Y=(Y(t))_{t\in[-\tau,T]}$ are solutions then $\mathbb{P}\big(\sup_{-\tau\leqslant s\leqslant t}\norm{Y(s)-X(s)}>0\big)=0$, for every $t\in[-\tau,T]$.
\end{definition}

The differentiability assumptions given in \cref{Section1} imply the following linear-growth and Lipschitz conditions. For each $j=0,\ldots,m$, there are constants $L_j^{(1)},L_j^{(2)}>0$ such that
\begin{align}
\norm{g_j(t,x_0,x_1,\ldots,x_K)}^2&\leqslant L_j^{(1)}\bigg(1+\sum_{k=0}^K\norm{x_k}^2\bigg),\label{LinearGrowthBound}\\
\norm{g_j(t,x_0,x_1,\ldots,x_K)-g_j(t,y_0,y_1,\ldots,y_K)}&\leqslant L_j^{(2)}\sum_{k=0}^K\norm{x_k-y_k},\label{LipschitzBound}
\end{align}
for all $(t,x_0,x_1,\ldots,x_K),(t,y_0,y_1,\ldots,y_K)\in[0,T]\times\mathbb{R}^d\times(\mathbb{R}^d)^K$, where $g_0=f$. The boundedness of $\sup_{t\in[-\tau,0]}\mathbb{E}[\norm{\phi(t)}^4]$ is required for the convergence order of some of our selected schemes (see the Magnus schemes in \cite{GriggsBurrageBurrage2025}). This boundedness of the fourth moment also implies $\sup_{t\in[-\tau,0]}\mathbb{E}[\norm{\phi(t)}^2]<\infty$, which, together with \eqref{LinearGrowthBound} and \eqref{LipschitzBound}, implies the existence of a unique solution for equation \eqref{Equation1}, due to the following result. This result is shown in the case of one delay, by K\"uchler and Platen \cite{KuchlerPlaten2000}, while the extension to $K$ delays uses similar reasoning.

\begin{theorem}\label{PLS}
Suppose that $f:[-\tau,T]\times\mathbb{R}^d\times(\mathbb{R}^d)^K\to\mathbb{R}^d$ is Lebesgue integrable, each $g_j:[-\tau,T]\times\mathbb{R}^d\times(\mathbb{R}^d)^K\to\mathbb{R}^d$ is square integrable, and also that each $f,g_j$ satisfies the linear-growth and Lipschitz bounds \eqref{LinearGrowthBound} and \eqref{LipschitzBound}. Additionally, suppose that $\phi$ is a stochastic process with $\int_{-\tau}^0\mathbb{E}[\norm{\phi(s)}^2]\,\mathrm{d}s<\infty$. Under these assumptions, there exists a pathwise-unique, strong solution $(X(t))_{t\in[-\tau,T]}$ of equation \eqref{Equation1}, satisfying $\sup_{t\in[-\tau,T]}\mathbb{E}[\norm{X(t)}^2]<\infty$.
\end{theorem}

\begin{remark}
\textcolor{black}{It is possible to loosen the conditions \eqref{LinearGrowthBound} and \eqref{LipschitzBound} while maintaining the existence and uniqueness of solutions for \eqref{Equation1}. For example, in the case of $K=1$ delay, Song, Hu, Gau, and Li \cite{SongHuGaoLi2022} show that the linear-growth bound may be replaced by a polynomial-growth bound. Alternatively, Przyby{\l}owicz, Wu, and Xie \cite{PrzybylowicsWuXie2024} show that the Lipschitz condition on the delayed argument $X(t-\tau_1)$ can be generalised to a H\"older condition where, for some $\alpha_j\in(0,1]$ and $\kappa>0$, $\norm{g_j(t,x_0,x_1)-g_j(t,y_0,y_1)}\leqslant \kappa\,(\norm{x_0-y_0}+\norm{x_1-y_1}^{\alpha_j})$, for each $j=0,\ldots,m$.}
\end{remark}


Equation \eqref{Equation1} is solved segmentally, between the multiples of the delay times. We list the multiples of the delays,
\[\begin{array}{ccccccc}
    \tau_1, & 2\tau_1, & 3\tau_1, & \cdots & \cdots & \cdots & N_1\tau_1, \\
    \tau_2, & 2\tau_2, & 3\tau_2, & \cdots & N_2\tau_2, & \; & \;  \\
    \vdots  &  \vdots  &  \vdots  & \cdots &    \ddots  & \; & \;  \\
    \tau_K, & 2\tau_K, & 3\tau_K, & \cdots & \cdots &  N_K\tau_K, & \;
  \end{array}\]
where $N_k=\inf\{N\in\mathbb{N}:N\geqslant T/\tau_k\}$ is the number of multiples of $\tau_k$ needed to reach $T$. We then sort these values into ascending order, letting this set be $\mathcal{T}^0=(\,\sigma_0^0,\sigma_1^0,\ldots,\sigma_{N_\mathrm{total}^0}^0\,)$, where $\sigma_0^0=0$ and $N_\mathrm{total}^0=N_1+\cdots+N_K$. Finally, we remove both redundant entries and also values exceeding $T$ from $\mathcal{T}^0$, labeling this final set of points $\mathcal{T}=(\sigma_0,\ldots,\sigma_{N_\mathrm{total}})$, where $\sigma_0=0$ and $\sigma_{N_\mathrm{total}}=T$. Due to the contributions by Bellman and Cooke \cite{BellmanCooke1963} to the study of ordinary delay-differential equations, we call each interval $[\sigma_b,\sigma_{b+1}]$ a \emph{Bellman interval}, for $b=0,\ldots,N_\mathrm{total}-1$. 

\begin{example}
Let $T=3$ and $K=4$ with delay values $\tau_1=1$, $\tau_2=2$, $\tau_3=2/3$, and $\tau_4=\pi/2$. In this case, $(N_1,N_2,N_3,N_4)=(3,2,5,2)$, and
\[\mathcal{T}^0=\left(\sigma_0^0,\sigma_1^0,\ldots,\sigma_{12}^0\right)=(0,\,2/3,\,1,\,4/3,\,\pi/2,\,2,\,2,\,6/3,\,8/3,\,3,\,\pi,\,10/3,\,4).\]
With redundancies removed, this list becomes
\[\mathcal{T}=(\sigma_0,\sigma_1,\ldots,\sigma_7)=(0,\,\tau_3,\,\tau_1,\,2\tau_3,\,\tau_4,\,2\tau_1,\,4\tau_3,\,3\tau_1)=(0,\,2/3,\,1,\,4/3,\,\pi/2,\,2,\,8/3,\,3).\]
\end{example}

On the interval $[\sigma_b,\sigma_{b+1}]$, the delayed solutions $X(t-\tau_k)$ are already determined, reducing \eqref{Equation1} to an SDE. That is, by writing $\phi_{\sigma_b}(t)=(X(t-\tau_1),\ldots,X(t-\tau_K))$ for $t\in[-\tau,\sigma_b]$, then \eqref{Equation1} on the present Bellman interval becomes
\begin{align}
\mathrm{d}X(t)&=[A_0X(t)+f(t,X(t),\phi_{\sigma_b}(t))]\,\mathrm{d}t+\sum_{j=1}^m[A_jX(t)+g_j(t,X(t),\phi_{\sigma_b}(t))]\,\mathrm{d}W_j(t),\quad t\in[\sigma_b,\sigma_{b+1}],\label{Equation2}
\end{align}
equipped with $X(\sigma_b)=\phi_{\sigma_b}(\sigma_b)$. Equation \eqref{Equation2} is readily solved numerically, by SDE schemes such as the Taylor and Magnus methods. In this work, we focus mainly on the Taylor SDE schemes, while also applying Magnus schemes to demonstrate extendability of our work to alternative schemes. The Taylor SDE schemes that we consider are those by Maruyama \cite{Maruyama1955} and Milstein \cite{Milstein1975}, which are well studied, such as in the text by Kloeden and Platen \cite{KloedenPlaten1999}. The Magnus SDE schemes are derived from exponential ODE schemes by Magnus \cite{Magnus1954}, extended to the homogeneous Stratonovich SDE setting by Burrage and Burrage \cite{BurrageBurrage1999}, the homogeneous It\^o SDE setting by Kamm, Pagliarani, and Pasucci \cite{KammPagliaraniPascucci2021}, and to inhomogeneous SDEs by Yang, Burrage, Komori, Burrage, and Ding \cite{YangBurrageKomoriBurrageDing2021}. The Magnus schemes that we consider are those of convergence orders $1/2$ and $1$, to complement the schemes by Maruyama and Milstein. More recently, the Taylor and Magnus schemes have been extended further, to SDDEs. The EM and Milstein SDDE schemes for \eqref{Equation1} are derived by K\"uchler and Platen \cite{KuchlerPlaten2000} by iteratively applying SDE schemes across the Bellman intervals. The EM method is the simplest numerical SDDE scheme, but its convergence may fail under non-global Lipschitz conditions, so truncated variants have been introduced by Guo, Mao, and Yue \cite{GanWang2014} to address this complication. However, the truncated EM scheme remains limited to OoC $1/2$, while the Milstein method achieves OoC $1$. Magnus schemes have also been extended to \eqref{Equation1}, by Griggs, Burrage, and Burrage \cite{GriggsBurrageBurrage2025}.

\subsection{Numerical schemes}\label{Section2.3}

In this paper, the numerical schemes that we consider are the Euler--Maruyama (EM), Milstein, Magnus--Euler--Maruyama (MEM), and Magnus--Milstein (MM) approximations for the solution of equation \eqref{Equation1}. Corresponding to these schemes is a time mesh $(t_n)_{n\in\mathcal{M}}$ satisfying $t_{n_1}<t_{n_2}$ whenever $n_1,n_2\in\mathcal{M}$ and $n_1<n_2$, where the index set $\mathcal{M}\subseteq\mathbb{Z}$ includes $\pm p_k$ for natural numbers $p_k$ chosen to satisfy $t_{p_k}=\tau_k$ and $t_{-p_k}=-\tau_k$, for $k=1,\ldots,K$. We also require $0,p,N\in\mathcal{M}$, where $t_0=0$, $t_N=T$, and $t_p=\tau$, with $p=\max_{k=1,\ldots,K}\{p_k\}$. That is, the set $\mathcal{M}$ (and corresponding time mesh $(t_n)_{n\in\mathcal{M}}$) may contain additional points, but at the minimum, we assume that it contains $0,\pm p_1,\ldots,\pm p_K,N$, and also all necessary points so that the schemes are well defined. We also find it convenient to write $(t_n)_{n\in\mathcal{M}}=(t_n)_{n=-p,\ldots,N}$, in sequence notation. A numerical scheme for the SDDE \eqref{Equation1} is a sequence $Y=(Y(t_n))_{n=-p,\ldots,N}$ intended to approximate the solution $X=(X(t))_{t\in[-\tau,T]}$. That is, $Y(t_n)\approx X(t_n)$, for each $n=-p,\ldots,N$. We abbreviate $Y(t_n)=Y_n$, for all $n=-p,\ldots,N$. We also set $Y_n=\phi(t_n)$ whenever $t_n\leqslant0$, for each numerical approximation of \eqref{Equation1}. 

\begin{definition}\label{DefinitionsofSchemes}
Suppose $Y=(Y_n)_{n=-p,\ldots,N}$ is a numerical scheme, approximating the solution $X=(X(t_n))_{t\in[-\tau,T]}$ of the SDDE \eqref{Equation1}, where $Y_n=\phi(t_n)$ whenever $t_n\leqslant0$. For each step $n=0,\ldots,N-1$, Wiener indices $i,j=1,\ldots,m$, and delay indices $k,l=1,\ldots,K$, the \emph{step size} is $h_n=t_{n+1}-t_n$, the \emph{delayed times} are $t_n^{\tau_k}=t_n-\tau_k$ and $t_n^{\tau_l,\tau_k}=t_n-\tau_l-\tau_k$, the \emph{delayed scheme values} are $Y_n^{\tau_k}=Y(t_n-\tau_k)$ and $Y_n^{\tau_l,\tau_k}=Y(t_n-\tau_l-\tau_k)$, the \emph{Wiener increment} is $\Delta W_j(t_n,t_{n+1})=W_j(t_{n+1})-W_j(t_n)$, the \emph{iterated stochastic integral} is the It\^o integral \eqref{Iij} and the \emph{delayed iterated stochastic integral} $I_{ij}^{\tau_k}(t_n,t_{n+1})$ is as defined by \eqref{Iijtauk}. We also define the integrals
\begin{equation}I_{i0}(t_n,t_{n+1})=\int_{t_n}^{t_{n+1}}\int_{t_n}^s\,\mathrm{d}W_i(u)\,\mathrm{d}s\quad\textrm{and}\quad I_{0j}(t_n,t_{n+1})=\int_{t_n}^{t_{n+1}}\int_{t_n}^s\,\mathrm{d}u\,\mathrm{d}W_j(s),\quad i\textrm{ and }j=1,\ldots,m.\label{Ii0andI0j}\end{equation}
The indicator function $\mathbb{I}(A)$ is given by
\[\mathbb{I}(A)=\begin{cases}1&\textrm{if }A\textrm{ is true,}\\
0&\textrm{otherwise,}\end{cases}\]
the \emph{commutator} (or \emph{Lie bracket}) of square matrices $A$ and $B$ is $[A,B]=AB-BA$, the (spatial) \emph{Jacobian} of $g_j$, at $(t,x,x_{\tau_1},\ldots,x_{\tau_K})$, is
\[\nabla_xg_j(t,x,x_{\tau_1},\ldots,x_{\tau_K})=\begin{pmatrix}
                            \frac{\partial g_{1,j}}{\partial x_1}(t,x,x_{\tau_1},\ldots,x_{\tau_K}) & \cdots & \frac{\partial g_{1,j}}{\partial x_d}(t,x,x_{\tau_1},\ldots,x_{\tau_K}) \\
                            \vdots & \ddots & \vdots \\
                            \frac{\partial g_{d,j}}{\partial x_1}(t,x,x_{\tau_1},\ldots,x_{\tau_K}) & \cdots & \frac{\partial g_{d,j}}{\partial x_d}(t,x,x_{\tau_1},\ldots,x_{\tau_K}) \\
                          \end{pmatrix}\]
and the \emph{delayed Jacobian} of each $g_j$, at $(t,x,x_{\tau_1},\ldots,x_{\tau_K})$, is
\[\nabla_{x_{\tau_k}}g_j(t,x,x_{\tau_1},\ldots,x_{\tau_K})=\begin{pmatrix}
                            \frac{\partial g_{1,j}}{\partial x_1^{\tau_k}}(t,x,x_{\tau_1},\ldots,x_{\tau_K}) & \cdots & \frac{\partial g_{1,j}}{\partial x_d^{\tau_k}}(t,x,x_{\tau_1},\ldots,x_{\tau_K}) \\
                            \vdots & \ddots & \vdots \\
                            \frac{\partial g_{d,j}}{\partial x_1^{\tau_k}}(t,x,x_{\tau_1},\ldots,x_{\tau_K}) & \cdots & \frac{\partial g_{d,j}}{\partial x_d^{\tau_k}}(t,x,x_{\tau_1},\ldots,x_{\tau_K}) \\
                          \end{pmatrix},\]
for each $j=1,\ldots,m$ and $k=1,\ldots,K$, where we write
\[(t,x,x_{\tau_1},\ldots,x_{\tau_K})=\left(t,(x_1,\ldots,x_d)^\intercal,(x_1^{\tau_1},\ldots,x_d^{\tau_1})^\intercal,\ldots,(x_1^{\tau_K},\ldots,x_d^{\tau_K})^\intercal\right)\in[0,T]\times\mathbb{R}^d\times(\mathbb{R}^d)^K.\]
Further to this, we denote $\tilde{f}(t,x,x_{\tau_1},\ldots,x_{\tau_K})=f(t,x,x_{\tau_1},\ldots,x_{\tau_K})-\sum_{j=1}^mA_jg_j(t,x,x_{\tau_1},\ldots,x_{\tau_K})$. The scheme $Y$ is:
\begin{enumerate}
  \item the \emph{Euler--Maruyama} (EM) approximation if
                \begin{equation}Y_{n+1}=Y_n+[A_0Y_n+f(t_n,Y_n,Y_n^{\tau_1},\ldots,Y_n^{\tau_K})]\,h_n+\sum_{j=1}^m[A_jY_n+g_j(t_n,Y_n,Y_n^{\tau_1},\ldots,Y_n^{\tau_K})]\,\Delta W_j(t_n,t_{n+1});\label{EMScheme}\end{equation}
  \item the \emph{Milstein} approximation if
                \begin{align}Y_{n+1}&=Y_n+[A_0Y_n+f(t_n,Y_n,Y_n^{\tau_1},\ldots,Y_n^{\tau_K})]\,h_n+\sum_{j=1}^m[A_jY_n+g_j(t_n,Y_n,Y_n^{\tau_1},\ldots,Y_n^{\tau_K})]\,\Delta W_j(t_n,t_{n+1})\nonumber\\
                &\;+\sum_{i=1}^m\sum_{j=1}^m\big[A_j+\nabla_xg_j(t_n,Y_n,Y_n^{\tau_1},\ldots,Y_n^{\tau_K})\big]\,\big[A_iY_n+g_i(t_n,Y_n,Y_n^{\tau_1},\ldots,Y_n^{\tau_K})\big]\,I_{ij}(t_n,t_{n+1})\label{MilsteinScheme}\\
                &\;+\sum_{k=1}^K\mathbb{I}(t_n\geqslant \tau_k)\sum_{i=1}^m\sum_{j=1}^m\nabla_{x_{\tau_k}}g_j(t_n,Y_n,Y_n^{\tau_1},\ldots,Y_n^{\tau_K})\,\left[A_iY_n^{\tau_k}+g_i(t_n^{\tau_k},Y_n^{\tau_k},Y_n^{\tau_1,\tau_k},\ldots,Y_n^{\tau_K,\tau_k})\right]\,I_{ij}^{\tau_k}(t_n,t_{n+1});\nonumber\end{align}
  \item the \emph{Magnus--Euler--Maruyama} (MEM) approximation if the matrix exponential
  \[M_n^{[1]}(t_n,t_{n+1})=\exp\bigg(\bigg(A_0-\frac{1}{2}\sum_{i=1}^mA_i^2\bigg)\,h_n+\sum_{j=1}^mA_j\,\Delta W_j(t_n,t_{n+1})\bigg)\]
                is used to compute
                \begin{equation}Y_{n+1}=M_n^{[1]}(t_n,t_{n+1})\,\bigg\{Y_n+\tilde{f}(t_n,Y_n,Y_n^{\tau_1},\ldots,Y_n^{\tau_K})\,h_n+\sum_{j=1}^mg_j(t_n,Y_n,Y_n^{\tau_1},\ldots,Y_n^{\tau_K})\,\Delta W_j(t_n,t_{n+1})\bigg\};\label{MEMScheme}\end{equation}
  \item the \emph{Magnus--Milstein} (MM) approximation if
                \[M_n^{[2]}(t_n,t_{n+1})=\exp\bigg(\bigg(A_0-\frac{1}{2}\sum_{i=1}^mA_i^2\bigg)\,h_n+\sum_{j=1}^mA_j\,\Delta W_j(t_n,t_{n+1})+\frac{1}{2}\sum_{i=0}^m\sum_{j=i+1}^m[A_i,A_j](I_{ji}(t_n,t_{n+1})-I_{ij}(t_n,t_{n+1}))\bigg)\]
                is used to compute
                \begin{align}Y_{n+1}&=M_n^{[2]}(t_n,t_{n+1})\,\Bigg\{Y_n+\tilde{f}(t_n,Y_n,Y_n^{\tau_1},\ldots,Y_n^{\tau_K})\,h_n+\sum_{j=1}^mg_j(t_n,Y_n,Y_n^{\tau_1},\ldots,Y_n^{\tau_K})\,\Delta W_j(t_n,t_{n+1})\label{MMScheme}\\
                &\hspace{-0.5cm}+\sum_{i=1}^m\sum_{j=1}^m\left\{\nabla_xg_j(t_n,Y_n,Y_n^{\tau_1},\ldots,Y_n^{\tau_K})\,\left[A_iY_n+g_i(t_n,Y_n,Y_n^{\tau_1},\ldots,Y_n^{\tau_K})\right]-A_ig_j(t_n,Y_n,Y_n^{\tau_1},\ldots,Y_n^{\tau_K})\right\}\,I_{ij}(t_n,t_{n+1})\nonumber\\
                &\;+\sum_{k=1}^K\mathbb{I}(t_n\geqslant \tau_k)\sum_{i=1}^m\sum_{j=1}^m\nabla_{x_{\tau_k}}g_j(t_n,Y_n,Y_n^{\tau_1},\ldots,Y_n^{\tau_K})\,\left[A_iY_n^{\tau_k}+g_i(t_n^{\tau_k},Y_n^{\tau_k},Y_n^{\tau_1,\tau_k},\ldots,Y_n^{\tau_K,\tau_k})\right]\,I_{ij}^{\tau_k}(t_n,t_{n+1})\Bigg\}\,,\nonumber\end{align}
\end{enumerate}
for each $n=0,\ldots,N-1$.
\end{definition}

\begin{definition}
Suppose $Y=(Y_n)_{n=-p,\ldots,N}$ is a numerical approximation for the solution $X=X(t)_{t\in[-\tau,T]}$ of equation \eqref{Equation1}, with fixed step size $h=h_n=t_{n+1}-t_n$, for each $n=-p,\ldots,N-1$. The scheme $Y$ has (strong) \emph{order of convergence} (OoC) $q$ (or \emph{convergence order} $q$) if there is a constant $c>0$ such that
\[\mathbb{E}\left[\norm{Y_n-X(t_n)}\right]\leqslant ch^q,\]
for all $n=-p,\ldots,N$. Similarly, $Y$ has \emph{mean-square order of convergence} (MS OoC) $q$ if there exists $c>0$ satisfying
\[\left(\mathbb{E}\left[\norm{Y_n-X(t_n)}^2\right]\right)^{1/2}\leqslant ch^q,\]
for all $n=-p,\ldots,N$.
\end{definition}


If a numerical scheme $Y$ has MS OoC $q$ then it can be shown (with the H\"older inequality for probability spaces) that $Y$ also has strong OoC $q$ (see \cite[page 5]{MilsteinTretyakov2021}). For this reason, it is often convenient to establish MS OoC, which then implies strong OoC. The EM and MEM schemes have OoC $1/2$, and if all relevant integrals can be expressed exactly, then the Milstein and MM schemes have OoC $1$ \cite{MilsteinTretyakov2021,GriggsBurrageBurrage2025}. It is evident from \cref{DefinitionsofSchemes} that the EM and MEM schemes are simpler to implement than the Milstein and MM schemes, at the cost of convergence order.

\subsection{Simulating second-order It\^o integrals}

\begin{figure}[ht]
\centering
\begin{subfigure}[b]{0.49\linewidth}
  \includegraphics[width=\linewidth, trim=2 2 2 2, clip]{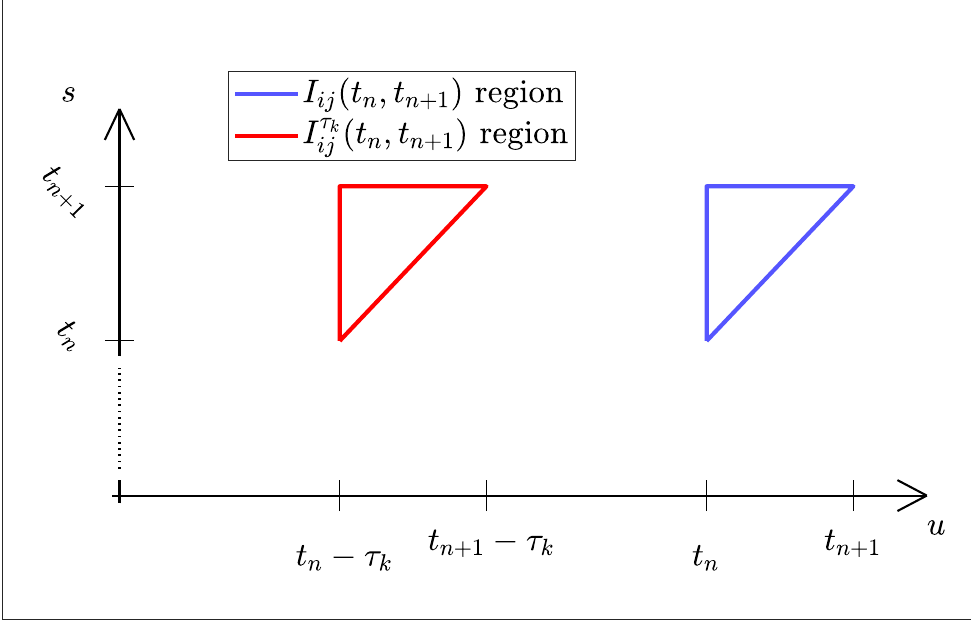}
  \caption{Integration regions of \eqref{Iij} and \eqref{Iijtauk}. These integrals can be approximated with \eqref{IijandIijtaukSimple}, using only the Wiener values at the corners of the integration regions.}
  \label{Figure2A}
\end{subfigure}
\hfill
\begin{subfigure}[b]{0.49\linewidth}
  \includegraphics[width=\linewidth, trim=2 2 2 2, clip]{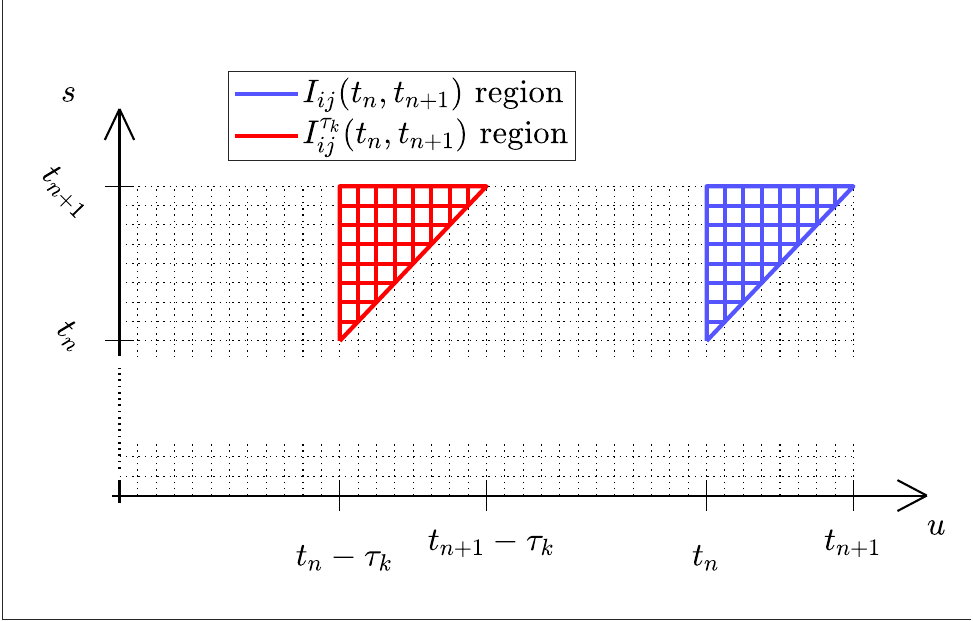}
  \caption{The integrals \eqref{Iij} and \eqref{Iijtauk} may also approximated with \eqref{IijandIijtaukRefinedApproximations}, using the Wiener values at time points in the refined mesh.}
  \label{Figure2B}
\end{subfigure}
\caption{Integration regions of $I_{ij}(t_n,t_{n+1})$ and $I_{ij}^{\tau_k}(t_n,t_{n+1})$ shown together on the same plane. To achieve strong order of convergence~1, the integrals are simulated using approximations \eqref{IijandIijtaukRefinedApproximations} based on the values of the Wiener paths evaluated on the finest time mesh. If the finest mesh uses a fixed step size $h_\mathrm{R}$ that divides $\tau_k$ then the delayed region also aligns with this mesh.}
\label{Figure2}
\end{figure}

In order to simulate the iterated integrals \eqref{Iij} and \eqref{Iijtauk}, we first observe that if we impose the bound
\begin{equation}\max_{n=0,\ldots,N}\{h_n\}<\min_{k=1,\ldots,K}\{\tau_k\}\label{Maxh}\end{equation}
on the largest possible step size for the numerical schemes, then even if $i=j$, the inner integral of \eqref{Iijtauk} has integration domain $[t_n-\tau_k,t_{n+1}-\tau_k]$ being disjoint from the outer domain $[t_n,t_{n+1}]$, for every $k=1,\ldots,K$. Therefore, if \eqref{Maxh} holds, then the independence of increments in the Wiener process ensures that the inner process $W_i^{\tau_k}(s)-W_i^{\tau_k}(t_n)=W_i(s-\tau_k)-W_i(t_n-\tau_k)$ and the outer process $W_j(s)$ are independent processes. This reduces the problem of simulating \eqref{Iijtauk} to the specific problem of simulating \eqref{Iij} under the condition \eqref{Maxh}. We further show this in \cref{Figure2}, where we show the integration regions of the integrals \eqref{Iij} and \eqref{Iijtauk}.

\begin{remark}
If some $h_n=t_{n+1}-t_n$ is greater than some delay $\tau_k$, then the integration regions overlap, so simulating the integrals without the condition \eqref{Maxh} requires further care. Fortunately, \eqref{Maxh} is not prohibitive as we are free to choose the step size for our schemes, and by selecting sufficiently small step sizes, the integration regions remain disjoint.
\end{remark}

Simulating \eqref{Iij} is complicated because the inner process of $\int_{t_n}^{t_{n+1}}\int_{t_n}^s\,\mathrm{d}W_i(u)\,\mathrm{d}W_j(s)=\int_{t_n}^{t_{n+1}}[W_i(s)-W_i(t_n)]\,\mathrm{d}W_j(s)$
 has unbounded variation, which renders quadrature-type approximations inadequate for achieving OoC $1$, but they suffice for OoC $1/2$. For example, we may simulate $\tilde{I}_{ij}\approx I_{ij}(t_n,t_{n+1})$ and $\tilde{I}_{ij}^{\tau_k}\approx I_{ij}^{\tau_k}(t_n,t_{n+1})$ with the approximations
\begin{equation}\tilde{I}_{ij}=\frac{\Delta W_i(t_n,t_{n+1})\,\Delta W_j(t_n,t_{n+1})}{2},\quad\tilde{I}_{ij}^{\tau_k}=\frac{\Delta W_i^{\tau_k}(t_n,t_{n+1})\,\Delta W_j(t_n,t_{n+1})}{2},\label{IijandIijtaukSimple}\end{equation}
which are \eqref{IijSimple} and \eqref{IijtaukSimple}, writing $\Delta W_i^{\tau_k}(s,t)=W_i^{\tau_k}(t)-W_i^{\tau_k}(s)=W_i(t-\tau_k)-W_i(s-\tau_k)$. The approximations \eqref{IijandIijtaukSimple} are easily inspired by \cref{Figure2A}, so we refer to these as \emph{simple} approximations of \eqref{Iij} and \eqref{Iijtauk}, and we call the scheme \eqref{MilsteinScheme} the \emph{simple Milstein scheme} when the integrals are simulated with \eqref{IijandIijtaukSimple}. That is, the simple Milstein scheme is the scheme by Hofmann and M\"uller-Gronbach \cite{HofmannMullerGronbach2006}, with OoC $1/2$.


In the case $i=j$, simulating $I_{jj}(t_n,t_{n+1})$ is simple, using the identity $I_{jj}(t_n,t_{n+1})=\big[(\Delta W_j(t_n,t_{n+1}))^2-h_n\big]/2$.
 When $i\neq j$, the Milstein scheme achieves OoC $1$ when the integration region is divided into subregions and we use the values of the Wiener processes at corresponding subpoints $t_n^{(l)}=t_n+lh_n/F_n$, for $l=0,\ldots,F_n$, where $F_n$ is a chosen number of subintervals within $[t_n,t_{n+1}]$. Specifically, we suppose a \emph{refined (time) mesh} is used on $[t_n,t_{n+1}]$, with \emph{refined step size} $h_n^\mathrm{R}>0$ chosen such that $F_n=h_n/h_n^\mathrm{R}\in\mathbb{N}$. From this assumption, the subpoints are of the form
\[t_n^{(l)}=t_n+lh_n^\mathrm{R},\quad l=0,\ldots,F_n.\]
In \cref{Figure2B}, we show the integration regions for $I_{ij}(t_n,t_{n+1})$ and $I_{ij}^{\tau_k}(t_n,t_{n+1})$ with constant step size $h_n=h$, displayed over (two-dimensional) refined time mesh, with constant refined step size $h^\mathrm{R}=h_n^\mathrm{R}$ chosen by $F_n=8$. For general $F_n$, the integral $I_{ij}(t_n,t_{n+1})$ can always be expressed over a refined time mesh, by
\begin{align}I_{ij}(t_n,t_{n+1})&=\sum_{l=0}^{F_n-1}\int_{t_n^{(l)}}^{t_n^{(l+1)}}\int_{t_n^{(l)}}^s\,\mathrm{d}W_i(u)\,\mathrm{d}W_j(s)+\sum_{l=0}^{F_n-1}\int_{t_n^{(l+1)}}^{t_n^{(F_n)}}\int_{t_n^{(l)}}^{t_n^{(l+1)}}\,\mathrm{d}W_i(u)\,\mathrm{d}W_j(s)\nonumber\\
&=\sum_{l=0}^{F_n-1}I_{ij}(t_n^{(l)},t_n^{(l+1)})+\sum_{l=0}^{F_n-1}\Delta W_i(t_n^{(l)},t_n^{(l+1)})\,\Delta W_j(t_n^{(l+1)},t_n^{(F_n)}).\label{IijDecomposition}\end{align}
Even using the decomposition \eqref{IijDecomposition}, the problem of calculating the integral $I_{ij}(t_n,t_{n+1})$ over a triangular region still involves the integrals $I_{ij}(t_n^{(l)},t_n^{(l+1)})$, which are also over triangular regions. However, by using this decomposition, it is possible to simulate the Milstein scheme with OoC $1$ (see \cite[Theorem 1.5.4]{MilsteinTretyakov2021}), by employing approximations $\tilde{I}_{ij}\approx I_{ij}(t_n,t_{n+1})$ such as the \emph{rectangle method}
\begin{equation*}\tilde{I}_{ij}=\sum_{l=0}^{F_n-1}\Delta W_i(t_n^{(l)},t_n^{(l+1)})\,\Delta W_j(t_n^{(l+1)},t_n^{(F_n)})\end{equation*}
and the \emph{trapezoid method}
\begin{equation}\tilde{I}_{ij}=\sum_{l=0}^{F_n-1}\frac{\Delta W_i(t_n^{(l)},t_n^{(l+1)})\,\Delta W_j(t_n^{(l)},t_n^{(l+1)})}{2}+\sum_{l=0}^{F_n-1}\Delta W_i(t_n^{(l)},t_n^{(l+1)})\,\Delta W_j(t_n^{(l+1)},t_n^{(F_n)})\label{IijTrapezoid}\end{equation}
by Milstein and Tretyakov \cite{Milstein1975,Milstein1995,MilsteinTretyakov2021}. Due to the use of a refined time mesh, we refer to the approximation \eqref{IijTrapezoid} as the \emph{refined} approximation for $I_{ij}(t_n,t_{n+1})$, while the \emph{refined Milstein} and \emph{refined MM} schemes are \eqref{MilsteinScheme} and \eqref{MMScheme}, respectively, implemented using \eqref{IijTrapezoid} to approximate $I_{ij}(t_n,t_{n+1})$.

\begin{remark}
The MM scheme \eqref{MMScheme} includes the integrals \eqref{Ii0andI0j}, which are simulated similarly to \eqref{Iij}, by defining $W_0(t)=t$.
\end{remark}

In this paper, we restrict ourselves to the simple and refined approximations of iterated stochastic integrals. However, alternative methods to simulate $I_{ij}(t_n^{(l)},t_n^{(l+1)})$ are also available. Kloeden, Platen, and Wright \cite{KloedenPlatenWright1992} expand this integral as a Fourier series, truncated after a number of terms. Wiktorsson \cite{Wiktorsson2001} improves this expansion by giving a precise expression for the remaining terms, rather than truncating them. Mrongowius and R\"{o}{\ss}ler \cite{MrongowiusRosler2021} improve this further, by providing a computationally faster expression for the remaining terms. Further expansions are given by Kuznetsov \cite{Kuznetsov2018,Kuznetsov2019}, using different Fourier coefficients.

Cao, Zhang, and Karniadakis \cite{CaoZhangKarniadakis2015} extend \eqref{IijTrapezoid}, as well as other Fourier-series representations of $I_{ij}(t_n^{(l)},t_n^{(l+1)})$, to the delayed analogue in the case of $K=1$ delay, and apply these to simulate the Milstein scheme \eqref{MilsteinScheme} with fixed scheme step size $h_n=h$ that divides the delay $\tau$. Extending this to $K$ delays, we introduce the \emph{refined} approximations
\begin{align}
\tilde{I}_{ij}      &= \sum_{l=0}^{F_n-1}\frac{\Delta W_i(t_n^{(l)},t_n^{(l+1)})\,\Delta W_j(t_n^{(l)},t_n^{(l+1)})}{2}+\sum_{l=0}^{F_n-1}\Delta W_i(t_n^{(l)},t_n^{(l+1)})\,\Delta W_j(t_n^{(l+1)},t_n^{(F_n)}),\nonumber\\
\tilde{I}_{ij}^{\tau_k}  &= \sum_{l=0}^{F_n-1}\frac{\Delta W_i^{\tau_k}(t_n^{(l)},t_n^{(l+1)})\,\Delta W_j(t_n^{(l)},t_n^{(l+1)})}{2}+\sum_{l=0}^{F_n-1}\Delta W_i^{\tau_k}(t_n^{(l)},t_n^{(l+1)})\,\Delta W_j(t_n^{(l+1)},t_n^{(F_n)}),\label{IijandIijtaukRefinedApproximations}
\end{align}
to approximate $\tilde{I}_{ij}\approx I_{ij}(t_n,t_{n+1})$ and $\tilde{I}_{ij}^{\tau_k}\approx I_{ij}^{\tau_k}(t_n,t_{n+1})$. By using these integral approximations, we are now equipped with the following collection of numerical approximations for equation \eqref{Equation1}:
\begin{enumerate}
  \item the EM scheme \eqref{EMScheme},
  \item the simple Milstein scheme \eqref{MilsteinScheme} using \eqref{IijandIijtaukSimple},
  \item the MEM scheme \eqref{MEMScheme},
  \item the simple MM scheme \eqref{MMScheme} using \eqref{IijandIijtaukSimple},
  \item the refined Milstein scheme \eqref{MilsteinScheme} using \eqref{IijandIijtaukRefinedApproximations}, and
  \item the refined MM scheme \eqref{MMScheme} using \eqref{IijandIijtaukRefinedApproximations}.
\end{enumerate}
The first four of these schemes have OoC $1/2$, while the last two have OoC $1$. These convergence orders may be proven using a result by Milstein \cite{Milstein1995}, which we present below (expressed for the delayed equation \eqref{Equation1}). Importantly, the refined schemes have the same convergence orders as their theoretical counterparts (which assume that we are able to simulate $I_{ij}(t_n,t_{n+1})$ and $I_{ij}^{\tau_k}(t_n,t_{n+1})$ exactly), so we use the refined schemes in place of the exact schemes in \cref{DefinitionsofSchemes}. We summarise these convergence orders in \cref{Table1}.


\begin{table}[H]
\centering
\begin{tabular}{|c|cccc|}
\hline
\multicolumn{5}{|c|}{\textbf{Convergence Orders for Different Integration Methods}} \\
\hline
&\multicolumn{4}{|c|}{Numerical Approximation} \\
Integration Method & EM & Milstein & MEM & MM \\
\hline
Exact Integrals & $1/2$ & $1$ & $1/2$ & $1$ \\
Refined Approximation & $1/2$ & $1$ & $1/2$ & $1$ \\
Simple Approximation & $1/2$ & $1/2$ & $1/2$ & $1/2$ \\
\hline
\end{tabular}
\caption{Orders of convergence for selected numerical schemes, using different methods to simulate the iterated stochastic integrals. Simple approximations are restricted to convergence order $1/2$, while refined approximations achieve the same convergence orders as the theoretical schemes, which assume exact simulation of all integrals.}
\label{Table1}
\end{table}

\begin{theorem}\label{GeneralTheoremofMilsteinforSDDEs}
Let $Y=(Y_n)_{n=-p,\ldots,N}$ be a numerical scheme with fixed step size $h=t_{n+1}-t_n$, approximating the solution $X=(X(t))_{t\in[-\tau,T]}$ of
equation \eqref{Equation1}, and suppose each delay $\tau_k=t_{p_k}$ for some mesh time $t_{p_k}$. If there are constants $M>0$, $q_2\geqslant 1/2$, and $q_1\geqslant q_2+1/2$, such that
\begin{align}
\norm{\mathbb{E}[Y_{n+1}|Y_n=x]-\mathbb{E}[X(t_{n+1})|X(t_n)=x]}&\leqslant M\left(1+\norm{x}^2\right)^{1/2}h^{q_1}\quad\textrm{and}\label{WeakBound}\\
\left(\mathbb{E}\left[\norm{Y_{n+1}-X(t_{n+1})}^2|Y_n=X(t_n)=x\right]\right)^{1/2}&\leqslant M\left(1+\norm{x}^2\right)^{1/2}h^{q_2},\label{MSPart}
\end{align}
for any $n=-p,\ldots,N-1$ and $x\in\mathbb{R}^d$, then
\begin{equation}\left(\mathbb{E}\left[\norm{Y_n-X(t_n)}^2|Y_0=X(0)\right]\right)^{1/2}\leqslant M\left(1+\mathbb{E}\left[\norm{X(0)}^2\right]\right)^{1/2}h^{q_2-1/2}\label{MSBoundSDDE}\end{equation}
holds for every $n=-p,\ldots,N$.
\end{theorem}

\subsection{Numerical demonstrations of convergence orders with divisible delays}

We now give examples demonstrating the convergence orders of our selected numerical schemes in the case of equation \eqref{Equation1} with $K=2$ divisible delays $\tau_1$ and $\tau_2$. In this setting, the delays have common factors, which we choose to be scheme step sizes. That is, there exists a step size $h>0$ along with natural numbers $p_1,p_2$ such that $\tau_1=h p_1$ and $\tau_2=h p_2$. We simulate the numerical schemes on the time mesh $(nh)_{n=-p,\ldots,T/h}$, while selecting a divisor $h^\mathrm{R}$ of $h$ and using the refined mesh $(nh^\mathrm{R})_{n=-p,\ldots,T/h^\mathrm{R}}$, on which we simulate a reference solution $X$, using the Milstein scheme, and the independent Wiener processes $W_1,\ldots,W_m$. We also write $p=\max\{p_1,p_2\}$, to complement $\tau=\max\{\tau_1,\tau_2\}$.

To show the convergence orders, we simulate each numerical approximation $Y=(Y_n)_{-p,\ldots,N}$, and then compute the mean-square error (MSE) $E_{t_n}=(\mathbb{E}\left\Vert Y_n-X(T)\right\Vert^2)^{1/2}$ between $Y_n=Y(t_n)$ and the solution value $X(t_n)$, at time $t_n$, using a large number $n_\mathrm{t}$ of trials. That is, we compute the error in the following steps. For each $i=1,\ldots,n_\mathrm{t}$:
\begin{enumerate}
  \item simulate the independent Wiener paths at the time points $nh^\mathrm{R}$, for $n=0,\ldots,T/h^\mathrm{R}$,
  \item simulate a reference solution $X^{(i)}=(X^{(i)}(t_n))_{n=-p,\ldots,T/h^\mathrm{R}}$ of \eqref{Equation1} with the Milstein scheme and fixed step size $h^\mathrm{R}$,
  \item simulate the numerical approximation $Y^{(i)}=(Y_n^{(i)})_{n=-p,\ldots,T/h}$ of \eqref{Equation1} with step size $h$,
\end{enumerate}
and then calculate
\begin{equation}
\mathrm{MSE}(t_n)=\left(\frac{1}{n_\mathrm{t}}\sum_{i=1}^{n_\mathrm{t}}\norm{Y^{(i)}_n-X^{(i)}(t_n)}^2\right)^{1/2}\approx E_{t_n}\label{MSE}
\end{equation}
as an approximation for $E_{t_n}$, where $Y_n^{(i)}=Y^{(i)}(t_n)$ is the numerical approximation at $t_n$ for the $i^\mathrm{th}$ trial. We calculate \eqref{MSE} for a range of step sizes $h$, to show how this error varies with $h$, on a $\log_2(h)$-$\log_{10}(\mathrm{MSE})$ error plot, for each numerical scheme. We also plot reference error graphs that show the behaviours of the theoretical convergence orders $1/2$ and $1$. To test the consistency of convergence behaviour across different delay configurations, we simulate these errors for three pairs of delays $(\tau_1,\tau_2)$.

\begin{example}\label{DivisibleExample}
We simulate the numerical schemes for equation \eqref{Equation1} with $m=2$ Wiener processes, spatial dimension $d=2$, terminal time $T=4$, and history process given by
\begin{equation}\phi(t)=\frac{1}{5}\begin{pmatrix} 4+t^2\sin(3\pi t) \\ 1+t^2\cos(2\pi t) \\ \end{pmatrix},\quad t\in[-\tau,0],\label{phit}\end{equation}
recalling that $\tau=\max\{\tau_1,\tau_2\}$ is the largest delay in equation \eqref{Equation1}. We select the matrices $A_0,A_1,A_2$ and functions $f,g_1,g_2$ given by
\begin{align}
A_0&=\begin{pmatrix}-0.1&0.03\\-0.2&-0.04\end{pmatrix},&f\left(t,\begin{pmatrix}x_1\\x_2\end{pmatrix},\begin{pmatrix}y_1\\y_2\end{pmatrix},\begin{pmatrix}z_1\\z_2\end{pmatrix}\right)&=\frac{1}{5}\begin{pmatrix}\sin(x_1)\\ \cos(x_2)\end{pmatrix},\nonumber\\
A_1&=\begin{pmatrix}0.05&0.04\\0.02&0.03\end{pmatrix},&g_1\left(t,\begin{pmatrix}x_1\\x_2\end{pmatrix},\begin{pmatrix}y_1\\y_2\end{pmatrix},\begin{pmatrix}z_1\\z_2\end{pmatrix}\right)&=\frac{1}{3}\begin{pmatrix}z_1-y_1\\y_2- z_2\end{pmatrix},\nonumber\\
A_2&=\begin{pmatrix}0.05&0.03\\0.04&0.01\end{pmatrix},&g_2\left(t,\begin{pmatrix}x_1\\x_2\end{pmatrix},\begin{pmatrix}y_1\\y_2\end{pmatrix},\begin{pmatrix}z_1\\z_2\end{pmatrix}\right)&=\frac{1}{10}\begin{pmatrix}\vphantom{\sum_k}\exp(-x_2^2)+\exp(-y_1^2)+\exp(-y_2^2)\\\vphantom{\sum_k}\exp(-x_1^2)+\exp(-z_1^2)+\exp(-z_2^2)\end{pmatrix}.\label{Example1Parameters}
\end{align}
The process \eqref{phit} is chosen to demonstrate non-constant history, as is typical in realistic applications. We show error graphs for three different pairs of delays to illustrate robustness (with respect to delay values) of the numerical schemes. The refined step size that we use to simulate the reference solution and Wiener processes is $h^\mathrm{R}=2^{-13}$. \Cref{ErrorGraphs1} shows the error graphs for the delay pairs $(\tau_1,\tau_2)=(1,1/2),\,(1,1/8),\,(1/4,1/8)$, with each graph showing the mean-square errors for the step sizes $h=\min\{\tau_1,\tau_2\},\ldots,2^{-10}$. These graphs are simulated with $n_\mathrm{t}=1,000$ trials.


\begin{figure}[ht]
    \centering
    \begin{subfigure}[t]{0.32\linewidth}
        \centering
        \includegraphics[width=\linewidth,trim=0 0 0 0,clip]{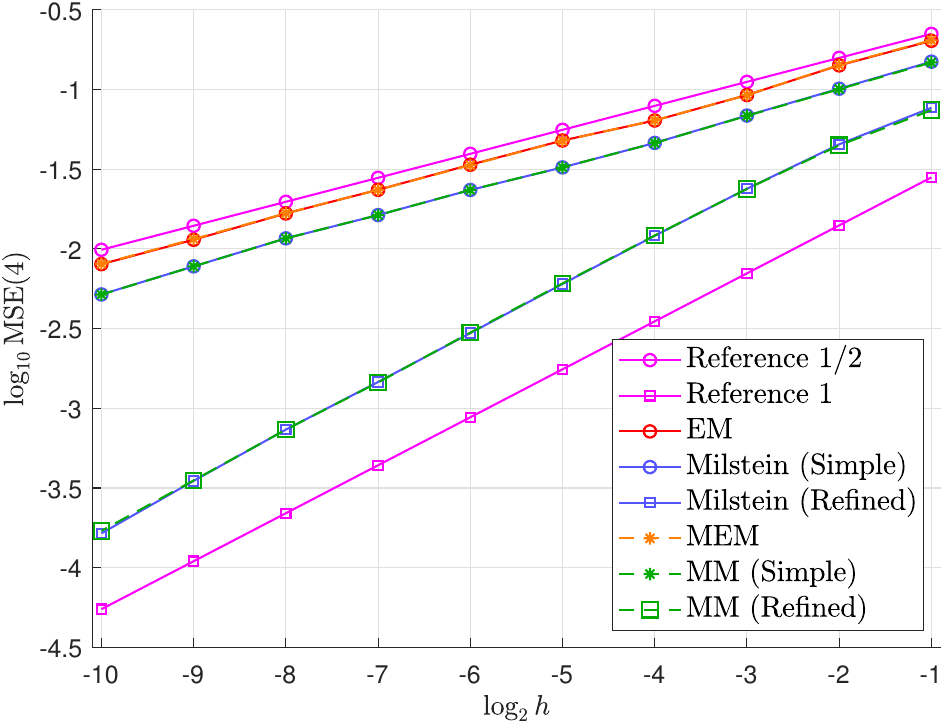}
        \caption{Error graphs with delays $\tau_1=1$ and $\tau_2=1/2$.}
        \label{Example1A}
    \end{subfigure}
    \hfill
    \begin{subfigure}[t]{0.32\linewidth}
        \centering
        \includegraphics[width=\linewidth,trim=0 0 0 0,clip]{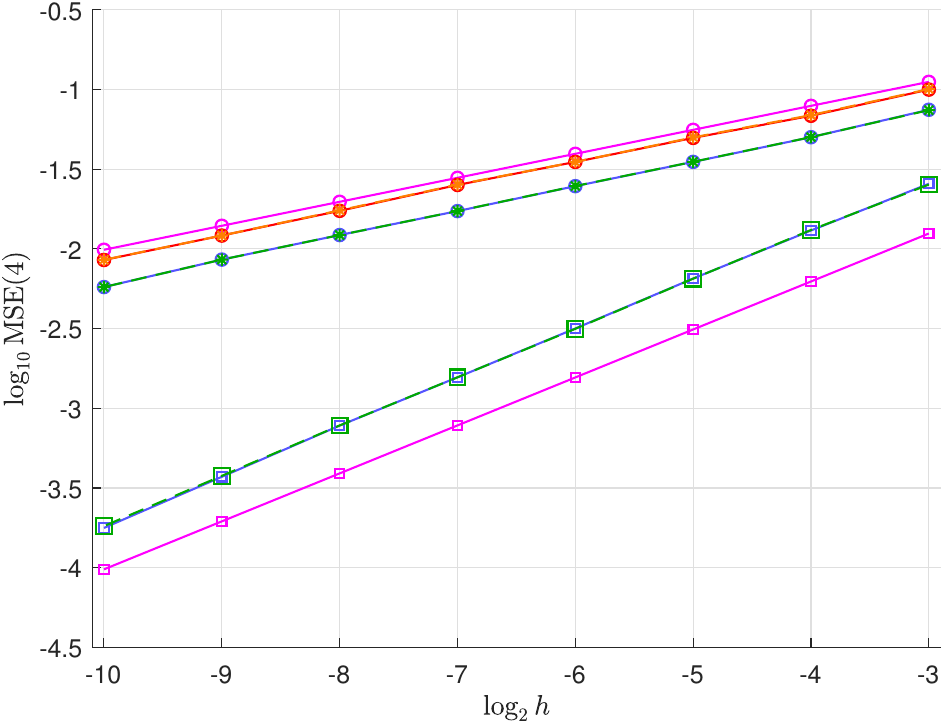}
        \caption{Error graphs with delays $\tau_1=1$ and $\tau_2=1/8$.}
        \label{Example1B}
    \end{subfigure}
    \hfill
    \begin{subfigure}[t]{0.32\linewidth}
        \centering
        \includegraphics[width=\linewidth,trim=0 0 0 0,clip]{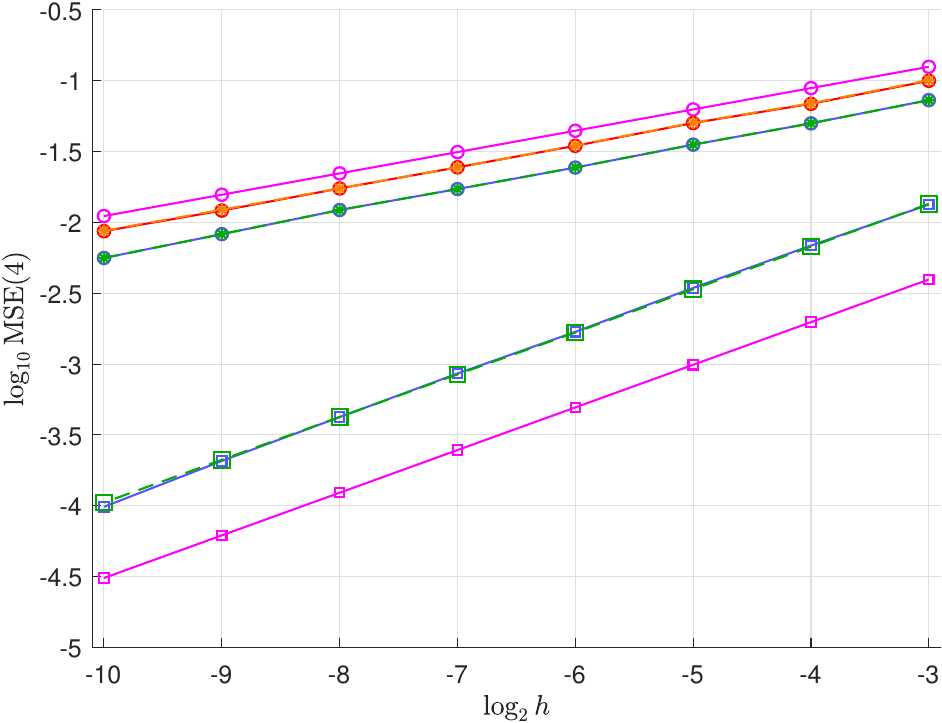}
        \caption{Error graphs with delays $\tau_1=1/4$ and $\tau_2=1/8$.}
        \label{Example1C}
    \end{subfigure}
    \caption{Mean-square error graphs for the numerical solutions of equation \eqref{Equation1} with the parameters given in \cref{DivisibleExample}, using three different pairs of delays $\tau_1$ and $\tau_2$. These numerical errors agree with the theoretical convergence orders shown in \cref{Table1}. Each error shown is evaluated at the terminal time $t_n=T=4$.}
    \label{ErrorGraphs1}
\end{figure}

\Cref{ErrorGraphs1} shows consistency in the convergence orders for the numerical schemes, across various delay values. The scheme errors are similar in \cref{Example1A,Example1B}, while \cref{Example1C} shows decreased errors in the schemes of OoC $1$, \eqref{MilsteinScheme} and \eqref{MMScheme}. This suggests that reducing the delay magnitudes may improve the numerical accuracy of SDDE schemes. Establishing this requires further investigation, but it may imply ramifications for applications of SDDE models.

\end{example}

\section{Fundamental complications for fixed-Step-size schemes achieving OoC 1 with indivisible delays}\label{Section3}

We now discuss the challenges of simulating numerical approximations with OoC $1$, for the SDDE \eqref{Equation1}. Theoretical existence of Milstein schemes for equation \eqref{Equation1} has been established by authors including Kloeden and Shardlow \cite{KloedenShardlow2012}, and Hu, Mohammed, and Yan \cite{HuMohammedYan2004}, in the case of general delays. However, practical implementation strategies remain underdeveloped, which we address in this paper. To begin this, we observe from \cref{DefinitionsofSchemes} that the schemes require simulation of delayed values such as $Y_n^{\tau_k}\approx X(t_n-\tau_k)$, as well as the integrals $I_{ij}(t_n,t_{n+1})$ and $I_{ij}^{\tau_k}(t_n,t_{n+1})$. We discuss the challenges of simulating these values in two parts: first, the delayed values; second, the integrals. For each of these challenges, we describe the simulation method in the case of divisible delays, before then considering the complications that arise when extending to indivisible delays.

\subsection{Continuous Milstein schemes}\label{ContinuousMilsteinSchemesSection}

Kloeden and Shardlow \cite{KloedenShardlow2012} give a continuous Milstein scheme $\overline{Y}$ for the SDDE
\begin{align}
\mathrm{d}X(t)&=\underline{f}(t,(X(s))_{s\in[t-\tau,t]})\,\mathrm{d}t+\sum_{j=1}^m\underline{g}_j(t,(X(s))_{s\in[t-\tau,t]})\,\mathrm{d}W_j(t),\quad t\in[0,T],\label{Equation3}\\
X(t)&=\underline{\phi}(t),\quad t\in[-\tau,0],\nonumber
\end{align}
where the history process $\underline{\phi}$ is continuous and $\underline{f},\underline{g}_1,\ldots,\underline{g}_m$
 satisfy certain differentiability conditions. Kloeden and Shardlow also prove that the continuous Milstein scheme has OoC $1$ in a continuous sense. Further to this, equation \eqref{Equation3} reduces to \eqref{Equation1} when the dependence (in each $\underline{f},\underline{g}_1,\ldots,\underline{g}_m$) on $(X(s))_{s\in[t-\tau,t]}$ is restricted to dependence only on the points $X(t),X(t-\tau_1),\ldots,X(t-\tau_K)$ (see \cite[Section 7]{KloedenShardlow2012}). For our purposes, we summarise these findings in the following theorem.

\begin{theorem}\label{TheoremKloedenShardlow}
Consider the SDDE
\begin{align*}
\mathrm{d}X(t)&=f(t,X(t),X(t-\tau_1),\ldots,X(t-\tau_K))\,\mathrm{d}t+\sum_{j=1}^mg_j(t,X(t),X(t-\tau_1),\ldots,X(t-\tau_K))\,\mathrm{d}W_j(t),\quad t\in[0,T],\\
X(t)&=\phi(t),\quad t\in[-\tau,0],
\end{align*}
with finitely many discrete delays $\tau_1,\ldots,\tau_K$, and coefficients $f$ and $g_j$ each satisfying \eqref{LinearGrowthBound}. Suppose the history process $\phi$ is continuous and the step size $h_n=t_{n+1}-t_n$ satisfies \eqref{Maxh}. On a time mesh $(t_n)_{t=-p,\ldots,N}$ including both $t_{p_k}=\tau_k$ and $t_{-p_k}=-\tau_k$, for some $p_k\in\mathbb{N}$, for each $k=1,\ldots,K$, suppose $Y=(Y_n)_{n=-p,\ldots,N}$ is the Milstein scheme defined by \eqref{MilsteinScheme}. The \emph{continuous Milstein scheme} $\overline{Y}=(\overline{Y}(t))_{t\in[-\tau,T]}$, defined for each $n=0,\ldots,N-1$ by
\begin{align}\overline{Y}(t)&=Y_n+f(t_n,Y_n,Y_n^{\tau_1},\ldots,Y_n^{\tau_K})\,(t-t_n)+\sum_{j=1}^mg_j(t_n,Y_n,Y_n^{\tau_1},\ldots,Y_n^{\tau_K})\,\Delta W_j(t_n,t)\nonumber\\
                &\;+\sum_{i=1}^m\sum_{j=1}^m\nabla_xg_j(t_n,Y_n,Y_n^{\tau_1},\ldots,Y_n^{\tau_K})\,g_i(t_n,Y_n,Y_n^{\tau_1},\ldots,Y_n^{\tau_K})\,I_{ij}(t_n,t)\label{ContinuousMilsteinScheme}\\
                &\;+\sum_{k=1}^K\mathbb{I}(t_n\geqslant \tau_k)\sum_{i=1}^m\sum_{j=1}^m\nabla_{x_{\tau_k}}g_j(t_n,Y_n,Y_n^{\tau_1},\ldots,Y_n^{\tau_K})\,g_i(t_n^{\tau_k},Y_n^{\tau_k},Y_n^{\tau_1,\tau_k},\ldots,Y_n^{\tau_K,\tau_k})\,I_{ij}^{\tau_k}(t_n,t),\quad t\in[t_n,t_{n+1}],\nonumber\end{align}
satisfies $\big(\mathbb{E}\big[\sup_{t\in[-\tau,T]}\big\Vert X(t)-\overline{Y}(t)\big\Vert^2\big]\big)^{1/2}\leqslant C\max_{n=-p,\ldots,N-1}\{t_{n+1}-t_n\}$, for constant $C>0$.
%
\end{theorem}


In order to simulate the continuous Milstein scheme \eqref{ContinuousMilsteinScheme} on the continuum $[0,T]$, we require a method to simulate the terms $I_{ij}(t_n,t)$ and $I_{ij}^{\tau_k}(t_n,t)$ at all times $t\in[0,T]$, which requires simulating a vast number of random variables. To avoid this, we instead simulate the integrals at a collection of discrete times between $t_n$ and $t_{n+1}$. That is, we simulate on a (refined) time mesh, and the theory by Kloeden and Platen proves that the OoC for such a Milstein scheme is $1$.

\subsection{Delayed scheme values}

\begin{figure}[ht]
    \centering
    \begin{subfigure}[t]{0.49\linewidth}
        \centering
        \includegraphics[width=\linewidth,trim=2 2 2 2,clip]{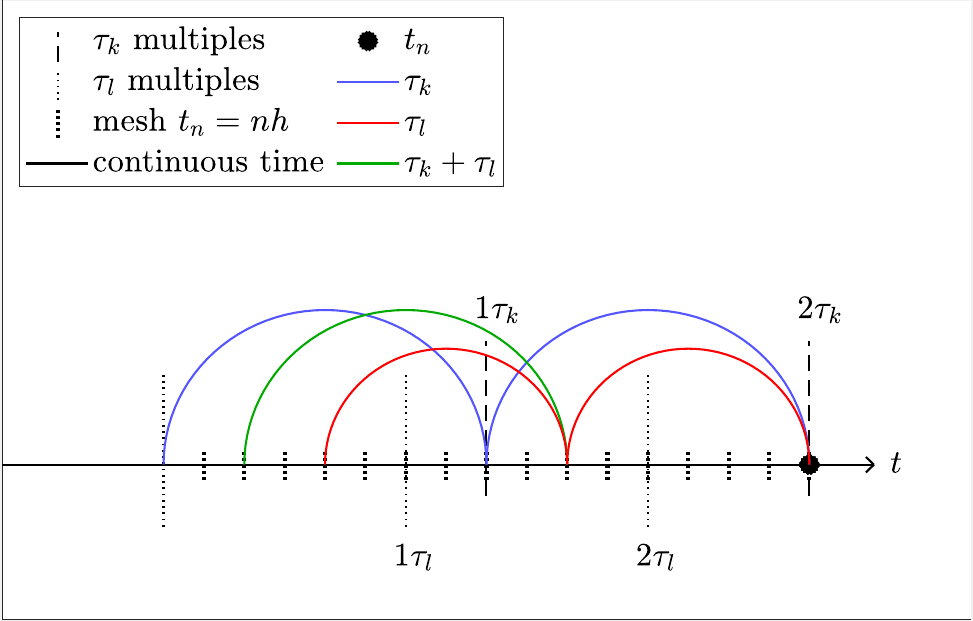}
        \caption{When the delays are all divisible by $h$, then every delayed time $t_n$, $t_n^{\tau_k}$, $t_n^{\tau_l,\tau_k}$ coincides with the time mesh, so that the corresponding scheme values $Y_n$, $Y_n^{\tau_k}$, $Y_n^{\tau_l,\tau_k}$ are already simulated.}
        \label{Figure1A}
    \end{subfigure}
    \hfill
    \begin{subfigure}[t]{0.49\linewidth}
        \centering
        \includegraphics[width=\linewidth,trim=2 2 2 2,clip]{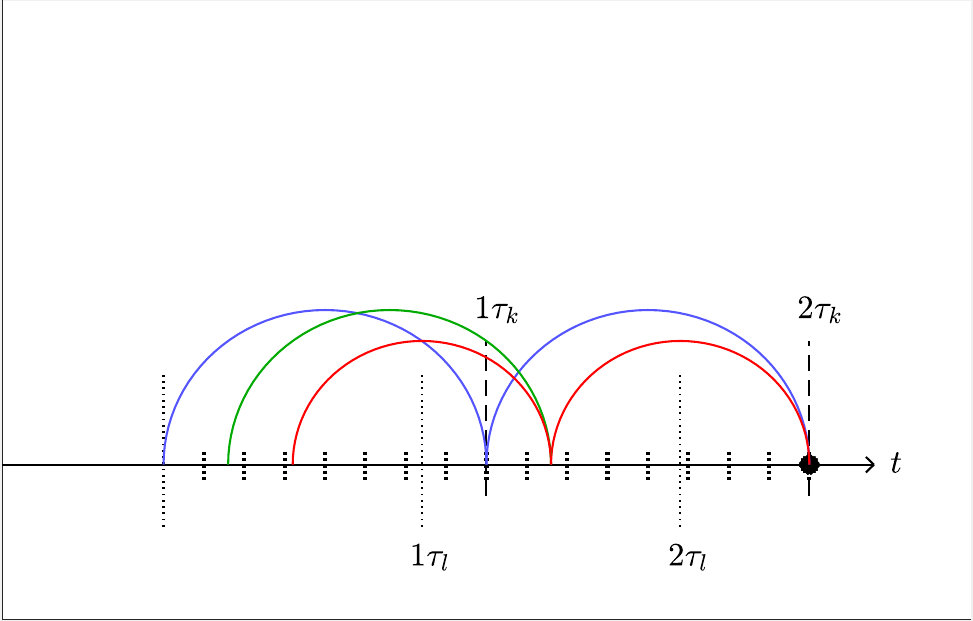}
        \caption{In the indivisible case, some delayed time $t_n^{\tau_l}$ is not on the time mesh, so that the scheme is not yet simulated at the delayed times $t_n^{\tau_l}$, $t_n^{\tau_l,\tau_k}$, $t_n^{\tau_l,\tau_l}$.}
        \label{Figure1B}
    \end{subfigure}
    \caption{Visual depiction of the (past and present) scheme values required to simulate a subsequent value, using the Milstein scheme. In order to compute $Y_{n+1}$ using \eqref{MilsteinScheme}, the scheme values $Y_n,\,Y_n^{\tau_k},\,Y_n^{\tau_l},\,Y_n^{\tau_k,\tau_k},\,Y_n^{\tau_l,\tau_k},\,Y_n^{\tau_l,\tau_l}$ are required, at respective times $t_n,\,t_n^{\tau_k},\,t_n^{\tau_l},\,t_n^{\tau_k,\tau_k},\,t_n^{\tau_l,\tau_k},\,t_n^{\tau_l,\tau_l}$. We visualise this on a time axis, starting at $t_n$, and showing delayed times $t_n-\tau_k$ and $t_n-2\tau_k$ (indicated with blue arcs), $t_n-\tau_l$ and $t_n-2\tau_l$ (red arcs), and the mixed delayed time $t_n-\tau_l-\tau_k$ (red followed by the green arc). In these figures, we assume that the current time $t_n$ coincides with $2\tau_k$, and also that $\tau_k$ aligns with the time mesh. The value $\tau_l$, however, no longer aligns with the time mesh, in \cref{Figure1B}.}
    \label{Figure1}
\end{figure}

The SDDE \eqref{Equation1} includes the delayed values $X(t-\tau_k)$, and similarly, each scheme requires the delayed values $Y_n^{\tau_k}$ in order to evaluate $Y_{n+1}$. If the delays are divisible by a common factor $h$ then we can apply a scheme using $h$ as the step size $h$. However, if the delays are indivisible, then there is no $h$ that divides every delay, meaning that a fixed-step-size scheme cannot simulate every delayed value $Y_n^{\tau_k}$. We show these situations in \cref{Figure1}, and to distinguish between these cases, we introduce the following definition.


\begin{definition}
We call the SDDE \eqref{Equation1} \emph{divisible} if there exists a common factor $h>0$ such that, for every $k=1,\ldots,K$, there is a natural number $p_k$ such that $hp_k=\tau_k$. Equation \eqref{Equation1} is \emph{indivisible} if the delays $\tau_1,\ldots,\tau_K$ have no common divisor.
\end{definition}

If the SDDE is divisible then a fixed step size $h=T/N$ may be used, where the time mesh $(t_n=nh)_{n=-p,\ldots,N}$ includes every delayed time $t_n-\tau_k=(n-p_k)\,h$ for $n=0,\ldots,N$. Therefore, when simulating $Y_{n+1}$, every delayed scheme value $Y_n^{\tau_k}=Y_{n-p_k}$ is already simulated. If the SDDE is indivisible, however, then we must either approximate some $Y_n^{\tau_k}$, or else we must loosen the requirement that $h$ be fixed.

\subsection{Simulating iterated stochastic integrals at arbitrary times}

\begin{figure}[ht]
\centering{\includegraphics[width=0.49\linewidth,trim=2 2 2 2,clip]{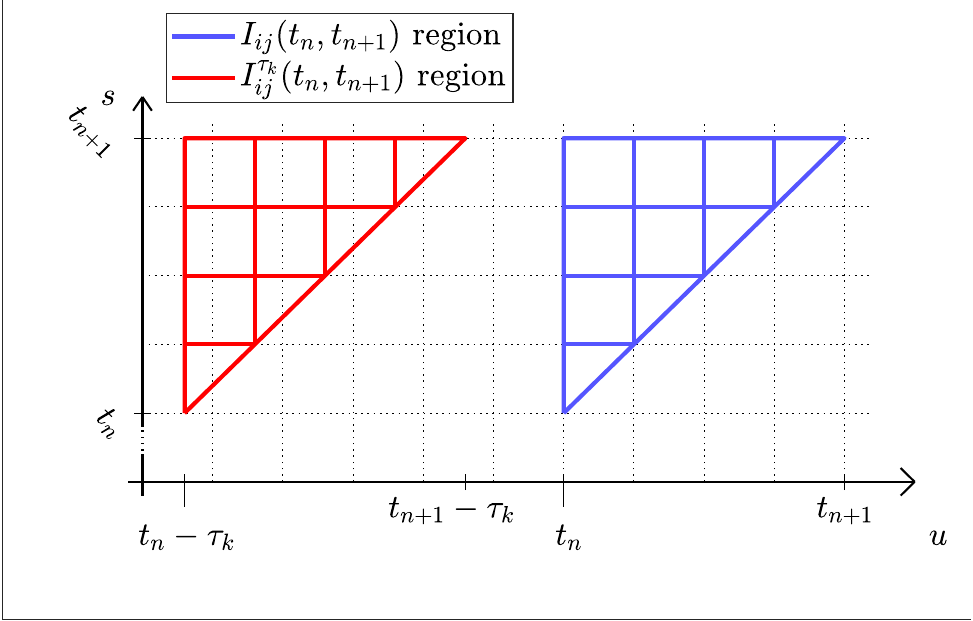}}
\caption{Integration regions of $I_{ij}(t_n,t_{n+1})$ and $I_{ij}^{\tau_k}(t_n,t_{n+1})$, shown over a refined time mesh of fixed step size $h^\mathrm{R}$, where the scheme step size $h=F_nh^\mathrm{R}$ (for natural number $F_n$) does not divide $\tau_k$. In this case, the approximation \eqref{IijandIijtaukRefinedApproximations} is not suited for $I_{ij}^{\tau_k}(t_n,t_{n+1})$, since the refined mesh does not coincide with the corners of the delayed integration region.}\label{Figure3}
\end{figure}

The continuous Milstein scheme requires simulation of integrals such as $I_{ij}(t_n,t)$ and $I_{ij}^{\tau_k}(t_n,t)$, for a continuum of values $t\in[t_n,t_{n+1}]$, which is generally complicated. 
 However, we are able to simulate these integrals at discrete points. In the case of a divisible SDDE, $I_{ij}(t_n,t_{n+1})$ and $I_{ij}^{\tau_k}(t_n,t_{n+1})$ may be simulated using \eqref{IijandIijtaukRefinedApproximations}. If instead the SDDE is indivisible, then a fixed-step-size mesh does not include some delayed time $t_n^{\tau_k}$, meaning that some Wiener value $W_i(t_n-\tau_k)$ is not yet simulated. This is also the case for the intermediate values $W_i(t_n^{(0)}-\tau_k),\ldots,W_i(t_n^{(F_n)}-\tau_k)$, which are not simulated if $\tau_k$ is not a multiple of $h$. We present this complication in \cref{Figure3}, showing that the delayed region no longer aligns with the underlying refined mesh.


\subsection{Discussion of solutions for these issues}

There are two practical approaches for addressing the issues discussed in this section, in the presence of indivisible delays:
\begin{enumerate}
  \item Use approximated delayed scheme values, such as computing $Y_n^{\tau_k}$ with linear interpolation between available scheme values. This approach is straightforward to implement and computationally efficient, but generally restricts the convergence order to $1/2$, which we see in the next section.
  \item Adapt the scheme to an augmented time mesh, constructed to include all required delay points. This allows us to evaluate delayed scheme values exactly, preserving the full order of convergence. This strategy is developed in \cref{Section5,Section6}, below.
\end{enumerate}





\begin{remark}
We have considered here complications that are ubiquitous across all strong-order-$1$ numerical approximations for \eqref{Equation1}, but individual schemes may include further complications. For example, a clear distinction between the Milstein and MM schemes is the presence of the matrix exponential $M_n^{[2]}(t_n,t_{n+1})$. Discussion of simulating \eqref{MMScheme} may be found in the derivation of the Magnus SDDE scheme in \cite{GriggsBurrageBurrage2025}, including methods for efficiently computing the matrix exponentials.
\end{remark}

\section{Method for strong OoC 1/2: Approximating delayed scheme values with linear interpolation}\label{Section4}

In this section, we establish a viable benchmark method to solve equation \eqref{Equation1} numerically when it is indivisible. We do this by maintaining a fixed step size while approximating the delayed scheme values using linear interpolation (LI) between the nearest simulated scheme values. This approach is computationally efficient and simple to implement, but it restricts the convergence order to $1/2$. That is, the LI method provides a simulation method for \eqref{Equation1}, offering efficient computation at the cost of convergence order.

\subsection{LI schemes}

\begin{figure}[ht]
\centering\includegraphics[width=0.49\linewidth,trim=4 4 4 4,clip]{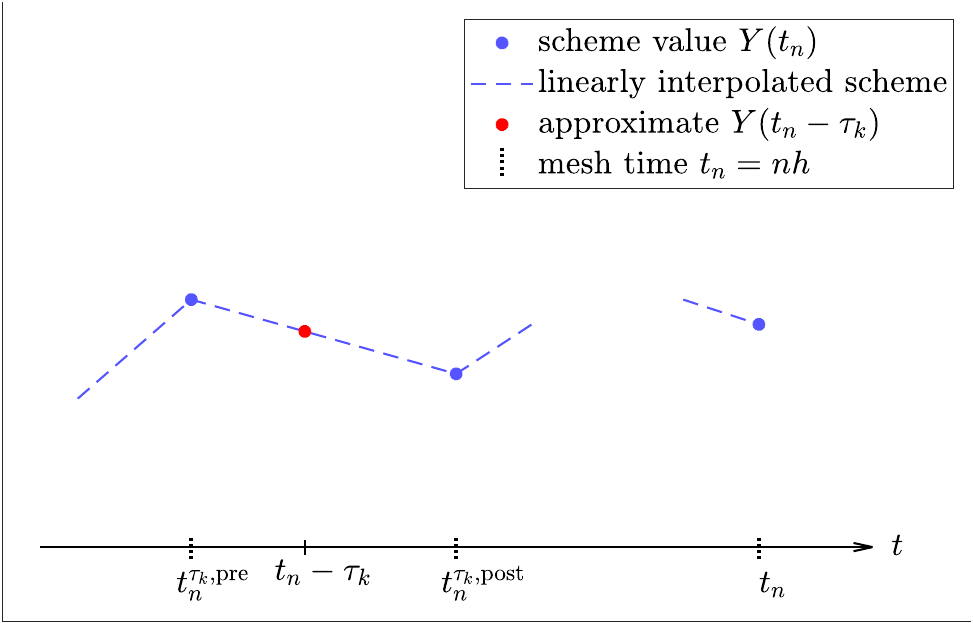}
\caption{Simulating the delayed scheme value $Y(t_n-\tau_k)$ using linear interpolation between the scheme values at the nearest mesh points. The nearest mesh points are those immediately before and after the delayed time $t_n-\tau_k$.}\label{Figure4}
\end{figure}

Suppose we are solving an indivisible equation \eqref{Equation1} with a numerical scheme of fixed step size $h$. The delayed time $t_n-\tau_k$ may land between mesh times, so we formulate the following general approach: we select greatest mesh time $t_n^{\tau_k,\mathrm{pre}}$ and least mesh time $t_n^{\tau_k,\mathrm{post}}$ satisfying $t_n^{\tau_k,\mathrm{pre}}\leqslant t_n-\tau_k$ and $t_n^{\tau_k,\mathrm{post}}>t_n-\tau_k$, and we then approximate the scheme $Y$ at the time $t_n-\tau_k$ using a linear interpolant between the selected mesh times. From prior simulation, we have already computed the corresponding scheme values $Y_n^{\tau_k,\mathrm{pre}}=Y(t_n^{\tau_k,\mathrm{pre}})$ and $Y_n^{\tau_k,\mathrm{post}}=Y(t_n^{\tau_k,\mathrm{post}})$, from which we form the linear approximation $\overline{Y}_n^{\tau_k}\approx X(t_n-\tau_k)$. We also define $\overline{Y}_n^{\tau_l,\tau_k}\approx X(t_n-\tau_l-\tau_k)$ similarly, by interpolating between the nearest mesh points around $t_n-\tau_l-\tau_k$. We provide visual representation for the LI construction of $\overline{Y}_n^{\tau_k}$ in \cref{Figure4}. 

\begin{definition}\label{DefinitionsofLISchemes}
Suppose $Y=(Y_n)_{n=-p,\ldots,N}$ is a numerical approximation of the solution $X=(X(t_n))_{t\in[-\tau,T]}$ for \eqref{Equation1}, with $Y_n=\phi(t_n)$ for $t_n\leqslant0$. Further to this, suppose $Y$ is a scheme of fixed step size $h_n=h$. For each step $n=0,\ldots,N-1$, we label the mesh times
\[t_n^{\tau_k,\mathrm{pre}}=\sup_{i=-p,\ldots,N}\left\{t_i:t_i\leqslant t_n-\tau_k\right\}\quad\textrm{and}\quad t_n^{\tau_k,\mathrm{post}}=\inf_{i=-p,\ldots,N}\left\{t_i:t_i>t_n-\tau_k\right\},\]
along with the scheme values $Y_n^{\tau_k,\mathrm{pre}}=Y(t_n^{\tau_k,\mathrm{pre}})$ and $Y_n^{\tau_k,\mathrm{post}}=Y(t_n^{\tau_k,\mathrm{post}})$, for each $k=1,\ldots,K$. Using these, we define the linear interpolant
\begin{equation}\overline{Y}_n^{\tau_k}=Y_n^{\tau_k,\mathrm{pre}}+\left(Y_n^{\tau_k,\mathrm{post}}-Y_n^{\tau_k,\mathrm{pre}}\right)\left(\frac{(t_n-\tau_k)-t_n^{\tau_k,\mathrm{pre}}}{h}\right),\quad n=0,\ldots,N-1,\quad k=1,\ldots,K.\label{Yntaukbar}\end{equation}
Similarly, the twice-delayed times are
\[t_n^{\tau_l,\tau_k,\mathrm{pre}}=\sup_{i=-p,\ldots,N}\left\{t_i:t_i\leqslant t_n-\tau_l-\tau_k\right\}\quad\textrm{and}\quad t_n^{\tau_l,\tau_k,\mathrm{post}}=\inf_{i=-p,\ldots,N}\left\{t_i:t_i>t_n-\tau_l-\tau_k\right\},\]
with twice-delayed scheme values $\overline{Y}_n^{\tau_l,\tau_k,\mathrm{pre}}=Y(t_n^{\tau_l,\tau_k,\mathrm{pre}})$ and $\overline{Y}_n^{\tau_l,\tau_k,\mathrm{post}}=Y(t_n^{\tau_l,\tau_k,\mathrm{post}})$, defined for $n$ when $t_n\geqslant\tau_k$. The twice-delayed linear interpolant is
\begin{equation}\overline{Y}_n^{\tau_l,\tau_k}=Y_n^{\tau_l,\tau_k,\mathrm{pre}}+\left(Y_n^{\tau_l,\tau_k,\mathrm{post}}-Y_n^{\tau_l,\tau_k,\mathrm{pre}}\right)\left(\frac{(t_n-\tau_l-\tau_k)-t_n^{\tau_l,\tau_k,\mathrm{pre}}}{h}\right),\quad n=0,\ldots,N-1,\quad l,k=1,\ldots,K.\label{Yntaultaukbar}\end{equation}
The \emph{linear interpolation} (LI) schemes are the schemes of \cref{DefinitionsofSchemes} constructed using the linear interpolants \eqref{Yntaukbar} and \eqref{Yntaultaukbar} in place of each $Y_n^{\tau_k}$ and $Y_n^{\tau_l,\tau_k}$, respectively. That is, $Y$ is:
\begin{enumerate}
  \item the \emph{Euler--Maruyama LI} (EM LI) scheme if
                \begin{equation*}Y_{n+1}=Y_n+[A_0Y_n+f(t_n,Y_n,\overline{Y}_n^{\tau_1},\ldots,\overline{Y}_n^{\tau_K})]\,h+\sum_{j=1}^m[A_jY_n+g_j(t_n,Y_n,\overline{Y}_n^{\tau_1},\ldots,\overline{Y}_n^{\tau_K})]\,\Delta W_j(t_n,t_{n+1});\end{equation*}
  \item the \emph{Milstein LI} (Milstein LI) scheme if
                \begin{align*}Y_{n+1}&=Y_n+[A_0Y_n+f(t_n,Y_n,\overline{Y}_n^{\tau_1},\ldots,\overline{Y}_n^{\tau_K})]\,h+\sum_{j=1}^m[A_jY_n+g_j(t_n,Y_n,\overline{Y}_n^{\tau_1},\ldots,\overline{Y}_n^{\tau_K})]\,\Delta W_j(t_n,t_{n+1})\nonumber\\
                &\;+\sum_{i=1}^m\sum_{j=1}^m\left[A_j+\nabla_xg_j(t_n,Y_n,\overline{Y}_n^{\tau_1},\ldots,\overline{Y}_n^{\tau_K})\right]\,\left[A_iY_n+g_i(t_n,Y_n,\overline{Y}_n^{\tau_1},\ldots,\overline{Y}_n^{\tau_K})\right]\,\tilde{I}_{ij}\\
                &\;+\sum_{k=1}^K\mathbb{I}(t_n\geqslant \tau_k)\sum_{i=1}^m\sum_{j=1}^m\nabla_{x_{\tau_k}}g_j(t_n,Y_n,\overline{Y}_n^{\tau_1},\ldots,\overline{Y}_n^{\tau_K})\,\left[A_i\overline{Y}_n^{\tau_k}+g_i(t_n^{\tau_k},\overline{Y}_n^{\tau_k},\overline{Y}_n^{\tau_1,\tau_k},\ldots,\overline{Y}_n^{\tau_K,\tau_k})\right]\,\tilde{I}_{ij}^{\tau_k};\nonumber\end{align*}
  \item the \emph{Magnus--Euler--Maruyama LI} (MEM LI) scheme if
  \[M_n^{[1]}(t_n,t_{n+1})=\exp\left(\left(A_0-\frac{1}{2}\sum_{i=1}^mA_i^2\right)\,h+\sum_{j=1}^mA_j\,\Delta W_j(t_n,t_{n+1})\right)\]
                and
                \begin{equation*}Y_{n+1}=M_n^{[1]}(t_n,t_{n+1})\,\left\{Y_n+\tilde{f}(t_n,Y_n,\overline{Y}_n^{\tau_1},\ldots,\overline{Y}_n^{\tau_K})\,h+\sum_{j=1}^mg_j(t_n,Y_n,\overline{Y}_n^{\tau_1},\ldots,\overline{Y}_n^{\tau_K})\,\Delta W_j(t_n,t_{n+1})\right\};\end{equation*}
  \item the \emph{Magnus--Milstein LI} (MM LI) scheme if the matrix exponential
                \[M_n^{[2]}(t_n,t_{n+1})=\exp\left(\left(A_0-\frac{1}{2}\sum_{i=1}^mA_i^2\right)\,h+\sum_{j=1}^mA_j\,\Delta W_j(t_n,t_{n+1})+\frac{1}{2}\sum_{i=0}^m\sum_{j=i+1}^m[A_i,A_j](\tilde{I}_{ji}-\tilde{I}_{ij})\right)\]
                is applied to compute
                \begin{align*}Y_{n+1}&=M_n^{[2]}(t_n,t_{n+1})\,\left\{\vphantom{\sum_{j=1}^m}\right.Y_n+\tilde{f}(t_n,Y_n,\overline{Y}_n^{\tau_1},\ldots,\overline{Y}_n^{\tau_K})\,h+\sum_{j=1}^mg_j(t_n,Y_n,\overline{Y}_n^{\tau_1},\ldots,\overline{Y}_n^{\tau_K})\,\Delta W_j(t_n,t_{n+1})\\
                &\hspace{-0.5cm}+\sum_{i=1}^m\sum_{j=1}^m\left\{\nabla_xg_j(t_n,Y_n,\overline{Y}_n^{\tau_1},\ldots,\overline{Y}_n^{\tau_K})\,\left[A_iY_n+g_i(t_n,Y_n,\overline{Y}_n^{\tau_1},\ldots,\overline{Y}_n^{\tau_K})\right]-A_ig_j(t_n,Y_n,\overline{Y}_n^{\tau_1},\ldots,\overline{Y}_n^{\tau_K})\right\}\,\tilde{I}_{ij}\nonumber\\
                &\;+\sum_{k=1}^K\mathbb{I}(t_n\geqslant \tau_k)\sum_{i=1}^m\sum_{j=1}^m\nabla_{x_{\tau_k}}g_j(t_n,Y_n,\overline{Y}_n^{\tau_1},\ldots,\overline{Y}_n^{\tau_K})\,\left[A_i\overline{Y}_n^{\tau_k}+g_i(t_n^{\tau_k},\overline{Y}_n^{\tau_k},\overline{Y}_n^{\tau_1,\tau_k},\ldots,\overline{Y}_n^{\tau_K,\tau_k})\right]\,\tilde{I}_{ij}^{\tau_k}\left.\vphantom{\sum_{j=1}^m}\right\}\,,\nonumber\end{align*}
\end{enumerate}
for each $n=0,\ldots,N-1$, where $\tilde{I}_{ij}$ and $\tilde{I}_{ij}^{\tau_k}$ are approximations of $I_{ij}(t_n,t_{n+1})$ and $I_{ij}^{\tau_k}(t_n,t_{n+1})$, respectively.
\end{definition}

\begin{remark}
If the SDDE \eqref{Equation1} is divisible, each $\tau_k=p_kh$ for some $p_k\in\mathbb{N}$, which implies that both $t_n^{\tau_k,\mathrm{pre}}=t_n^{\tau_k}$ and $t_n^{\tau_k,\mathrm{post}}=t_n^{\tau_k}+h$, while $\overline{Y}_n^{\tau_k}=Y_{n-p_k}=Y_n^{\tau_k}$. Similarly, each $\overline{Y}_n^{\tau_l,\tau_k}=Y_{n-p_l-p_k}=Y_n^{\tau_l,\tau_k}$. Therefore the LI schemes reduce to the usual schemes from \cref{DefinitionsofSchemes}, in the case of divisible SDDEs.
\end{remark}

\subsection{Numerical examples}

We now demonstrate the convergence orders of the numerical schemes using linear interpolation to simulate delayed scheme values. To show these convergence orders, we replicate \cref{DivisibleExample}, but with different pairs of delay values. We simulate the LI schemes and again calculate the error \eqref{MSE}, comparing each scheme with a reference solution $X$. We simulate the reference solution using the Milstein scheme on the \emph{augmented mesh}, which we develop in \cref{Section5,Section6}. We also simulate the Wiener processes and the stochastic integrals using this augmented mesh.

\begin{example}\label{LIExample}

We simulate the EM, simple Milstein, refined Milstein, MEM, simple MM, and refined MM schemes using linear interpolants to give delayed scheme values. The parameters $m,d,T$ are as in \cref{DivisibleExample}, with $\phi$ given by \eqref{phit}, and matrices $A_0,A_1,A_2$ and functions $f,g_1,g_2$ given by \eqref{Example1Parameters}. We again simulate the schemes for a range of step sizes $h=2^{-1},\ldots,2^{-10}$, with the errors shown at time $t=T=4$, on a $\log_2(h)$-$\log_{10}(\mathrm{MSE})$ error plot, produced with $n_\mathrm{t}=2,500$ trials. We show the error graphs for these numerical approximations, in \cref{ErrorGraphs1LI}, for the pairs of delays $(\tau_1,\tau_2)=(1,1/2),\,(1,\pi/4),\,(\exp(2)/5,\pi/4)$. To distinguish these schemes from those in \cref{DefinitionsofSchemes}, we label them with LI in \cref{ErrorGraphs1LI}.

\begin{figure}[ht]
    \centering
    \begin{subfigure}[t]{0.32\linewidth}
        \centering
        \includegraphics[width=\linewidth,trim=0 0 0 0,clip]{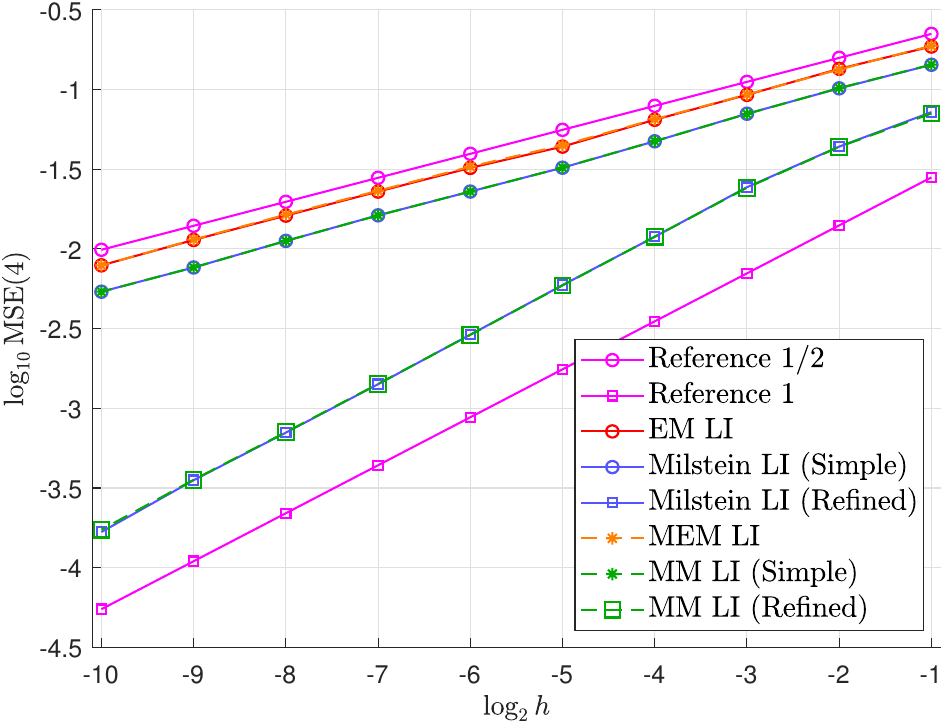}
        \caption{Error graphs with delays $\tau_1=1$ and $\tau_2=1/2$.}
        \label{Example1LIA}
    \end{subfigure}
    \hfill
    \begin{subfigure}[t]{0.32\linewidth}
        \centering
        \includegraphics[width=\linewidth,trim=0 0 0 0,clip]{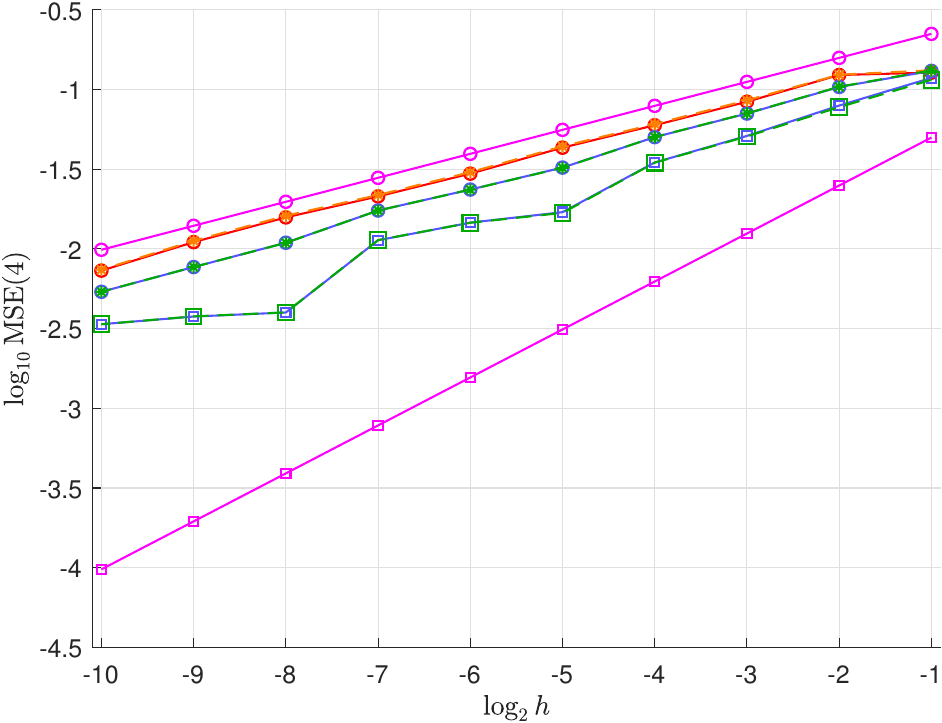}
        \caption{Error graphs with $\tau_1=1$ and $\tau_2=\pi/4$.}
        \label{Example1LIB}
    \end{subfigure}
    \hfill
    \begin{subfigure}[t]{0.32\linewidth}
        \centering
        \includegraphics[width=\linewidth,trim=0 0 0 0,clip]{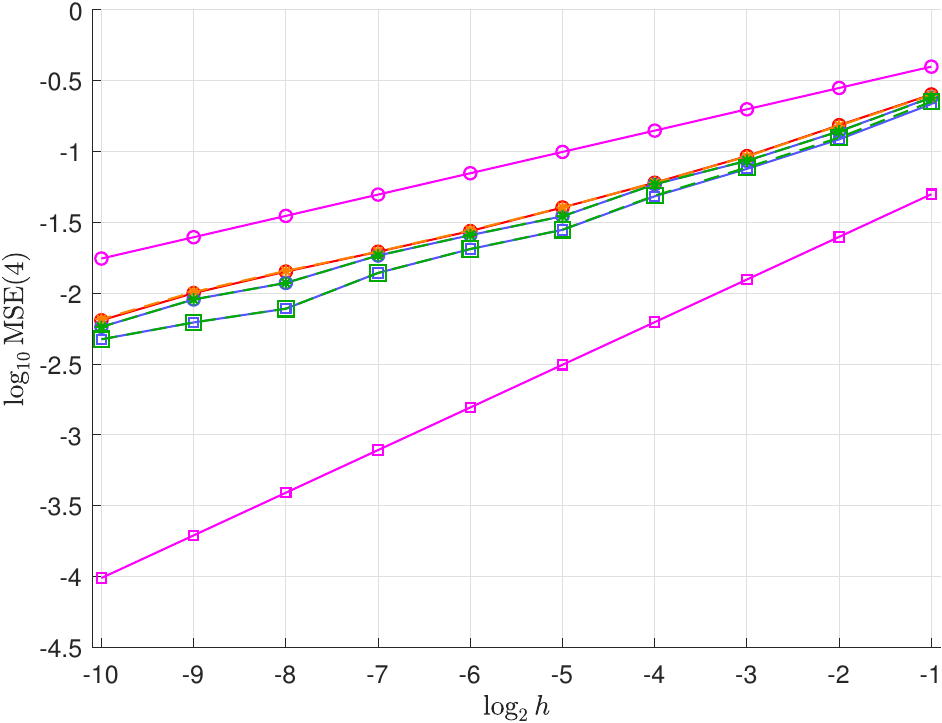}
        \caption{Error graphs with $\tau_1=\exp(2)/5$ and $\tau_2=\pi/4$.}
        \label{Example1LIC}
    \end{subfigure}
    \caption{Mean-square error graphs for the fixed-step LI numerical solutions of equation \eqref{Equation1} with the parameters given in \cref{DivisibleExample}, using three different pairs of delays $\tau_1$ and $\tau_2$. When the delays are not on the scheme mesh, linear interpolation is used to simulate the delayed scheme values.}
    \label{ErrorGraphs1LI}
\end{figure}

In these examples, we have simulated the reference solution (using \emph{initial refined step size} $h_\mathrm{initial}^\mathrm{R}=2^{-13}$), Wiener processes, and importantly the integral approximations $\tilde{I}_{ij}\approx I_{ij}(t_n,t_{n+1})$ and $\tilde{I}_{ij}^{\tau_k}\approx I_{ij}^{\tau_k}(t_n,t_{n+1})$ as accurately as possible. We later see in \cref{Section7} that OoC $1$ is achieved by the refined schemes constructed in this accurate manner. However, we see in \cref{ErrorGraphs1LI} that although the error graphs for the LI schemes display OoC $1$ (in the case of the refined schemes) when the delays are divisible (\cref{Example1LIA}), the schemes display OoC $1/2$ when the delays are indivisible (\cref{Example1LIB,Example1LIC}). We therefore conclude that the convergence order of the refined Milstein (and MM) scheme is broken when approximating delayed scheme values with linear interpolation, if the delays are indivisible.


\end{example}

The LI schemes provide computationally efficient methods to solve equation \eqref{Equation1} numerically. However, we see from \cref{LIExample} that in the case of an indivisible SDDE, each of the numerical schemes is restricted to OoC $1/2$, so if OoC $1/2$ is the goal in this example, then there is little reason to use anything more complicated than the EM LI scheme. More generally, however, there may still be reasons to use other schemes. For example, Griggs, Burrage, and Burrage \cite{GriggsBurrageBurrage2025} find that when numerically solving a delayed stochastic heat equation (DSHE) with fine spatial discretisation, the MEM scheme is numerically stable where the EM scheme is not, for larger temporal step sizes. By combining this observation with our conclusions above, the MEM LI scheme is better suited than the EM LI scheme, for numerically solving a DSHE with indivisible delays.

\subsection{Discussion of alternative methods using approximated delayed values}


We have presented the LI schemes to provide a benchmark method for solving equation \eqref{Equation1} numerically. However, future research may investigate alternative numerical methods. We conclude this section by discussing some of these potential techniques.
\begin{enumerate}
  \item \textbf{Left-hand values:} Rather than applying linear interpolation, it is possible to define instead $Y(t_n-\tau_k)=Y_n^{\tau_k,\mathrm{pre}}$, using the previous simulated value before the time $t_n-\tau_k$. This provides minor computational savings (avoiding the computations \eqref{Yntaukbar} and \eqref{Yntaultaukbar}). In general, we expect this approach to only produce OoC $1/2$ since it involves no more work than the LI method.
  \item \textbf{Backwards-forwards SDEs:} Rather than using a linear segment between $Y(t_n^{\tau_k,\mathrm{pre}})$ and $Y(t_n^{\tau_k,\mathrm{post}})$ in \cref{Figure4}, we may instead simulate $Y(t_n-\tau_k)$ as the solution at $t=t_n-\tau_k$ of the backwards-forwards SDE
      \begin{align}
      \mathrm{d}Y(t)&=[A_0Y(t)+f(t,Y(t),Y(t-\tau_1),\ldots,Y(t-\tau_K))]\,\mathrm{d}t\nonumber\\
      &\quad+\sum_{j=1}^m[A_jY(t)+g_j(t,Y(t),Y(t-\tau_1),\ldots,Y(t-\tau_K))]\,\mathrm{d}W_j(t),\quad t\in[t_n^{\tau_k,\mathrm{pre}},t_n^{\tau_k,\mathrm{post}}],\label{BFSDE}
      \end{align}
      equipped with predetermined values $Y(t_n^{\tau_k,\mathrm{pre}})$ and $Y(t_n^{\tau_k,\mathrm{post}})$. Ma and Yong \cite[page 254]{MaYong2007} show that it is possible to solve a backwards-forwards SDE numerically with OoC $1/2$. However, to our knowledge, it is not known if there are higher-order numerical solutions for equation \eqref{BFSDE}.
  \item \textbf{Collocation techniques:} Extending further from the LI schemes, other continuous extensions may be used to define $Y(t_n-\tau_k)$, such as polynomial collocation methods, with work by Conte et al.\ \cite{Conte2021} or by Fodor, Sykora, and Bachrathy \cite{FodorSykoraBachrathy2023}.
\end{enumerate}

\section{Method for strong OoC 1: Enhanced-mesh discretisation}\label{Section5}

In the remainder of this paper, we focus on simulating the (refined) Milstein scheme \eqref{MilsteinScheme}, with the other schemes simply provided as examples in implementations. We have already established simulation methods for the Milstein scheme when equation \eqref{Equation1} is divisible, so we now turn our attention towards indivisible SDDEs. We derive an implementation approach that realizes the theoretical convergence order given by \cref{TheoremKloedenShardlow}. We recall from our above work that the Milstein scheme achieving OoC 1 requires both an enhanced scheme mesh and also accurate simulations of the integrals \eqref{Iij} and \eqref{Iijtauk} (and also accurate simulations of \eqref{Ii0andI0j} in the case of the Magnus--Milstein scheme \eqref{MMScheme}). We address the mesh enhancement in this section, while detailing simulation of the integrals in \cref{Section6}.

We construct a scheme mesh that includes all time points required to simulate the Milstein scheme \eqref{MilsteinScheme}. We begin this construction by defining the set of non-delayed times at which the scheme must be evaluated, which we call the set of \emph{observation times}, labelled $T^\mathrm{O}$. For this work, we suppose that this set is
\begin{equation}T^\mathrm{O}=\{nh_\mathrm{initial}\}_{n=0,\ldots,T/h_\mathrm{initial}}\cup\{i_k\tau_k:i_k=1,\ldots,N_k,\,k=1,\ldots,K\},\label{TO}\end{equation}
for a fixed \emph{initial} step size $h_\mathrm{initial}>0$, and recalling that $N_k=\inf\{N\in\mathbb{N}:N\geqslant T/\tau_k\}$ is the number of multiples of $\tau_k$ required to pass $T$. With $h_\mathrm{initial}$ chosen sufficiently small (to satisfy \eqref{Maxh}) as a divisor of $T$, the set $T^\mathrm{O}$ given by \eqref{TO} includes a mesh of fixed step size combined with the delay multiples. We include the delay multiples $i_k\tau_k$ in $T^\mathrm{O}$ to ensure the indicator $\mathbb{I}(t_n\geqslant\tau_k)$ in \eqref{MilsteinScheme} is always evaluated unambiguously, and to capture the Bellman interval boundaries necessary for evaluation.

\begin{remark}
We have freedom to include additional observation times beyond those given in \eqref{TO}. For example, if we require simulation at time $t=\pi$ then we include $\pi$ in $T^\mathrm{O}$. We also may loosen the condition that $h_\mathrm{initial}$ divides $T$, but must include $T$ in the set $T^\mathrm{O}$.
\end{remark}

The set $T^\mathrm{O}$ contains times at which we seek numerical solutions, but the schemes require simulation at other, prior times. From the Milstein scheme \eqref{MilsteinScheme}, in order to simulate at a time $t_{n+1}$, we require prior simulation at times $t_n$, $t_n-\tau_k$, and also $t_n-\tau_l-\tau_k$, for each $l,k=1,\ldots,K$. However, in order to simulate the scheme at these points, the scheme requires prior simulation at $t_n-\tau_o-\tau_l-\tau_k$, for $o,l,k=1,\ldots,K$, and so on. Therefore the scheme requires simulation at all times in the set
\begin{equation*}T^\mathrm{C}=\{t_n-i_1\tau_1-i_2\tau_2-\cdots-i_K\tau_K:t_n\in T^\mathrm{O},\,i_k=0,\ldots,N_k,\,k=1,\ldots,K\}.\end{equation*}
We call $T^\mathrm{C}$ the (unordered) \emph{complete mesh}. We do not require simulation outside the time interval $[0,T]$, so we remove any values from $T^\mathrm{C}$ that are less than $0$. We also remove duplicate elements from $T^\mathrm{C}$, and then order the remaining elements in ascending order. We call the resulting ordered set the \emph{augmented (time) mesh}. We summarise this procedure in \cref{AlgorithmAugmentedMesh}.

\begin{algorithm}[H]
  \SetAlgoLined
  \LinesNumbered
  \KwData{SDDE \eqref{Equation1}, delay times $\tau_1,\ldots,\tau_K$, terminal time $T$, initial step size $h_\mathrm{initial}>0$}
  set observation times $T^\mathrm{O}=\{nh_\mathrm{initial}:n=0,1,\ldots,T/h_\mathrm{initial}\}\cup\{i\tau_k:i=1,\ldots,N_k,\,k=1,\ldots,K\}$\;
  set complete list $T^\mathrm{C}=\{t_n-i_1\tau_1-\cdots-i_K\tau_K:t_n\in T^\mathrm{O},\;i_k=0,\ldots,N_k,\,k=1,\ldots,K\}$\;
  remove from $T^\mathrm{C}$ any times less than $0$\;
  remove duplicate entries from $T^\mathrm{C}$\;
  sort $T^\mathrm{C}$ into ascending order\;
  \caption{Construction of the augmented time mesh for the Milstein scheme \eqref{MilsteinScheme}.}\label{AlgorithmAugmentedMesh}
\end{algorithm}

\begin{remark}\label{RemarkonGeneralisabilityforHigherOrderSchemes}
Higher-order numerical SDDE schemes include simulation at times such as $t_n-\tau_o-\tau_l-\tau_k$, or more generally, $t_n-\sum_{k=1}^Ki_k\tau_k$ for $i_k=0,1,\ldots$, similarly to the Milstein scheme \eqref{MilsteinScheme} including simulation at times $t_n-\tau_l-\tau_k$. The augmented mesh constructed by \cref{AlgorithmAugmentedMesh} includes all such points, so may also be used for these schemes.
\end{remark}

In practical implementation of \cref{AlgorithmAugmentedMesh}, caution is required to avoid excessively large arrays, due to line 2. If this step is performed iteratively across the delays (by first defining $T^\mathrm{C}=\{t_n-i_1\tau_1:t_n\in T^\mathrm{O},\,i_1=0,\ldots,N_1\}$ and then updating $T^\mathrm{C}=\{t_n-i_k\tau_k:t_n\in T^C,\,i_k=0,\ldots,N_k\}$ for each $k=2,\ldots,K$ separately), then $T^\mathrm{C}$ can become so large that computers are unable to store all values due to memory restrictions. This memory issue can be avoided by removing duplicate values within each iteration. For example, the pseudocode \cref{AlgorithmOptimalMatlab} may be used to produce the augmented time mesh (lines 2--5 in \cref{AlgorithmAugmentedMesh}), in the case of $K=4$ delays. Although line 8 of \cref{AlgorithmOptimalMatlab} could be deferred until after the loops without altering the final result, applying it within each iteration reduces the size of the vector, which assists with reducing memory usage during construction.

\begin{algorithm}[H]
  \SetAlgoLined
  \LinesNumbered
  initialise the set $T^\mathrm{C}=T^\mathrm{O}$\;
  \For{$i_1=0,\ldots,N_1$}
       {\For{$i_2=0,\ldots,N_2$}
       {\For{$i_3=0,\ldots,N_3$}
       {\For{$i_4=0,\ldots,N_4$}
       {$T^\mathrm{C}=[T^\mathrm{C},T^\mathrm{C}-i_1\tau_1-i_2\tau_2-i_3\tau_3-i_4\tau_4]$\;
       remove elements from $T^\mathrm{C}$ less than $0$\;
       remove duplicate values from $T^\mathrm{C}$\;
       sort $T^\mathrm{C}$ into ascending order\;
            }
            }
            }
            }
  \caption{Snippet of pseudocode used for memory-saving construction of augmented mesh $T^\mathrm{C}$, in the case of $K=4$ delays.}\label{AlgorithmOptimalMatlab}
\end{algorithm}

\begin{example}
In this example, we enumerate the size of the augmented mesh constructed by \cref{AlgorithmAugmentedMesh}, for a range of initial step sizes. We fix $T=1$ and consider up to $K=4$ delays. Shown in \cref{Table2} are the sizes of the augmented mesh $\norm{T^\mathrm{C}}$, corresponding to the initial step sizes $h_\mathrm{initial}=2^{-2},2^{-4},2^{-6},2^{-8},2^{-10}$, for each listed collection of delays $(\tau_1,\tau_2,\tau_3,\tau_4)$. Six collections of delays are chosen, to demonstrate the following cases.
\begin{enumerate}
  \item All divisible delays, showing that the augmented mesh agrees with the initial mesh $(nh_\mathrm{initial})_{n=0,\ldots,T/h_\mathrm{initial}}$ for divisible SDDEs.
  \item Divisible delays again, where each delay remains on the initial mesh.
  \item One indivisible delay, which is large compared with $T$.
  \item Two large indivisible delays.
  \item Three large indivisible delays.
  \item Four indivisible delays, including small (compared with $T$) delays.
\end{enumerate}

\begin{table}[H]
\centering
  \begin{tabular}{|c|c|c|c|c|c|c|c|c|}
    \hline
    \multicolumn{9}{|c|}{\textbf{\vphantom{$\displaystyle{\int}$}Size of Augmented Mesh $\norm{T^\mathrm{C}}$}} \\
    \hline
    \multicolumn{4}{|c|}{\vphantom{$\displaystyle{\sum}$}Delays} & \multicolumn{5}{|c|}{\vphantom{$\displaystyle{\sum}$}Initial Step Size $h_\mathrm{initial}$} \\
    $\tau_1$ & $\tau_2$ & $\tau_3$ & $\tau_4$ & $2^{-2}$ & $2^{-4}$ & $2^{-6}$ & $2^{-8}$ & $2^{-10}$ \\
    \hline
    \vspace{-10pt}&&&&&&&&\\
    \hline
    $1$ & $1$ & $1$ & $1$                                   & $5$ & $17$ & $65$ & $257$ & $1025$ \\
    \hline
    $1/4$ & $1$ & $1$ & $1$                                 & $5$ & $17$ & $65$ & $257$ & $1025$ \\
    \hline
    $1/4$ & $\pi/4$ & $1$ & $1$                             & $10$ & $25$ & $83$ & $316$ & $1249$ \\
    \hline
    $1/4$ & $\pi/4$ & $1/\sqrt{6}$ & $1$                    & $25$ & $51$ & $154$ & $572$ & $2240$ \\
    \hline
    $1/4$ & $\pi/4$ & $1/\sqrt{6}$ & $1/\sqrt{3}$            & $35$ & $66$ & $190$ & $695$ & $2709$ \\
    \hline
    $1/10$ & $\pi/10$ & $1/\sqrt{10}$ & $\exp(-2)/2$            & $4344$ & $6688$ & $16,505$ & $55,519$ & $211,734$ \\
    \hline
  \end{tabular}
  \caption{Number of points $\norm{T^\mathrm{C}}$ in the augmented mesh, until terminal time $T=1$, constructed by \cref{AlgorithmAugmentedMesh}, given delays $\tau_1,\tau_2,\tau_3,\tau_4$ and initial step size $h_\mathrm{initial}$. The time points in the augmented mesh are of the form $t_n-\sum_{k=1}^4i_k\tau_k\geqslant0$, for integers $i_k\geqslant0$, corresponding to initial mesh points $t_n=nh_\mathrm{initial}$ with $n=0,\ldots,T/h_\mathrm{initial}$.}\label{Table2}
\end{table}

\cref{AlgorithmAugmentedMesh} generalises the initial mesh, since the latter agrees with the augmented mesh in the case of divisible delays, as demonstrated in the first two rows of \cref{Table2}. In the third row, the only indivisible delay $\pi/4\approx0.79$ is large compared with $T=1$, so the added points $t_n-\pi/4$ are only nonnegative for a few of the initial mesh points $t_n$. This also occurs in the fourth row, with the additional large delay $1/\sqrt{6}\approx0.41$, as well as in the fifth row, with large delay $1/\sqrt{3}\approx0.58$. In the final row, however, the delays $(1/10,\pi/10,1/\sqrt{10},\exp(-2)/2)\approx(0.1,0.31,0.32,0.07)$ include smaller delays, resulting in significant increase in the number of added points 
 to the augmented mesh.

\end{example}

From the enumerations shown in \cref{Table2}, it is evident that the practicality of the augmented mesh is determined by the ratios of the delays compared with the terminal time. When these ratios $\tau_k/T$ are large, the augmented mesh remains a viable simulation method. However, as the delays become small compared with $T$, the number of mesh points may become prohibitively large. \textcolor{black}{This follows from our construction, since the size of the set \eqref{TO} is bounded by $\norm{T^\mathrm{O}}\leqslant 1+T/h_\mathrm{initial}+\sum_{k=1}^KN_k$, and
 for each observation time ($t_n\in T^\mathrm{O}$), the augmented mesh also contains at most $\sum_{k=1}^K(N_k+1)$ additional values (of the form $t_n-\sum_{k=1}^K\sum_{i_k=0}^{N_k}i_k\tau_k$), which gives the upper bound
\begin{equation}\norm{T^\mathrm{C}}\leqslant\left(1+T/h_\mathrm{initial}+\sum_{k=1}^KN_k\right)\sum_{k=1}^K(N_k+1).\label{UpperBoundonTc}\end{equation}
Bounds that are tighter than \eqref{UpperBoundonTc} may be found, for example by removing redundancies and negative values (as in \cref{AlgorithmAugmentedMesh}), but in the general case of indivisible delays, such bounds remain quadratic in $N_1,\ldots,N_K$. This is seen in \cref{Table2}, where $\norm{T^\mathrm{C}}$ approximately increases squarely as each $N_k=\lceil T/\tau_k\rceil$ increases.}


We have now constructed an augmented time mesh to simulate numerical approximations for equation \eqref{Equation1}. In the case of the Milstein scheme \eqref{MilsteinScheme}, we also require simulation of the integrals \eqref{Iij} and \eqref{Iijtauk}, on this mesh. We address simulation of these integrals in \cref{Section6}, but before concluding this section, we include a remark useful for comparing the accuracy of numerical schemes built on augmented meshes.

\begin{remark}\label{RemarkonAugmentedMeshesBeingRefinements}
Suppose we are applying \cref{AlgorithmAugmentedMesh} with an initial step size $h_{\mathrm{initial}}^{(i)}$, to produce a complete mesh $T^{\mathrm{C},(i)}$. If we repeat this for a sequence of initial step sizes $h_{\mathrm{initial}}^{(1)},h_{\mathrm{initial}}^{(2)},\ldots$, where each $h_{\mathrm{initial}}^{(i+1)}$ divides $h_{\mathrm{initial}}^{(i)}$, then it follows that each $T^{C,(i)}$ is a subset of $T^{C,(i+1)}$. In particular, halving the initial step size produces a refinement of the augmented mesh.
\end{remark}

\section{Integration simulation on the augmented mesh}\label{Section6}

We extend the trapezoidal approximations \eqref{IijandIijtaukRefinedApproximations} to indivisible SDDEs, using the augmented mesh. In the case of divisible SDDEs, selecting a refined step size $h^\mathrm{R}$ that divides the scheme step size $h$ ensures that the refined mesh $(nh^\mathrm{R})_{n=0,\ldots,T/h^\mathrm{R}}$ contains the scheme mesh $(nh)_{n=0,\ldots,T/h}$ (where we are disregarding the negative-time mesh points). From this, simulating the Wiener processes on the refined mesh allows us to implement the approximations \eqref{IijandIijtaukRefinedApproximations}, which is required for the refined numerical schemes of OoC $1$. In order to extend this construction to indivisible SDDEs, we apply \cref{RemarkonAugmentedMeshesBeingRefinements}, and select an initial refined step size $h_\mathrm{initial}^\mathrm{R}$ that divides the initial step size $h_\mathrm{initial}$.

\begin{definition}
Suppose we are simulating a numerical solution for an indivisible SDDE, by applying the augmented mesh constructed with \cref{AlgorithmAugmentedMesh}, using an initial (fixed) step size $h_\mathrm{initial}$. If $h_\mathrm{initial}^\mathrm{R}>0$ divides $h_\mathrm{initial}$ then the \emph{augmented refined time mesh} (ARTM) is the augmented mesh (constructed with \cref{AlgorithmAugmentedMesh}) applied using $h_\mathrm{initial}^\mathrm{R}$ as the initial step size.
\end{definition}

\begin{figure}[ht]
    \centering
    \begin{subfigure}[t]{0.49\linewidth}
        \centering
        \includegraphics[width=\linewidth,trim=2 2 2 2,clip]{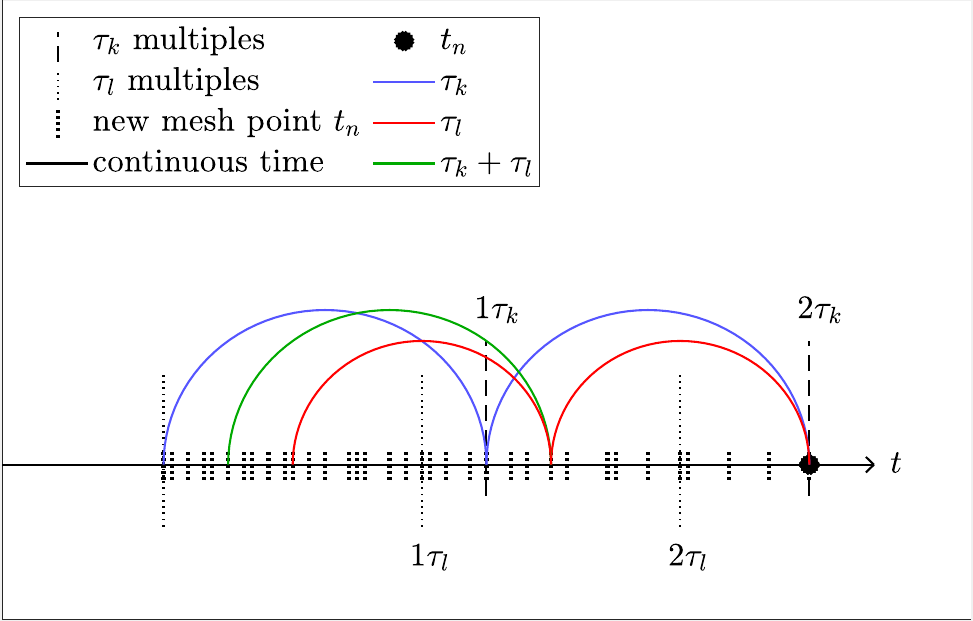}
        \caption{The augmented mesh includes all delayed times needed to compute $Y_{n+1}$.}
        \label{Figure1C}
    \end{subfigure}
    \hfill
    \begin{subfigure}[t]{0.49\linewidth}
        \centering
        \includegraphics[width=\linewidth,trim=2 2 2 2,clip]{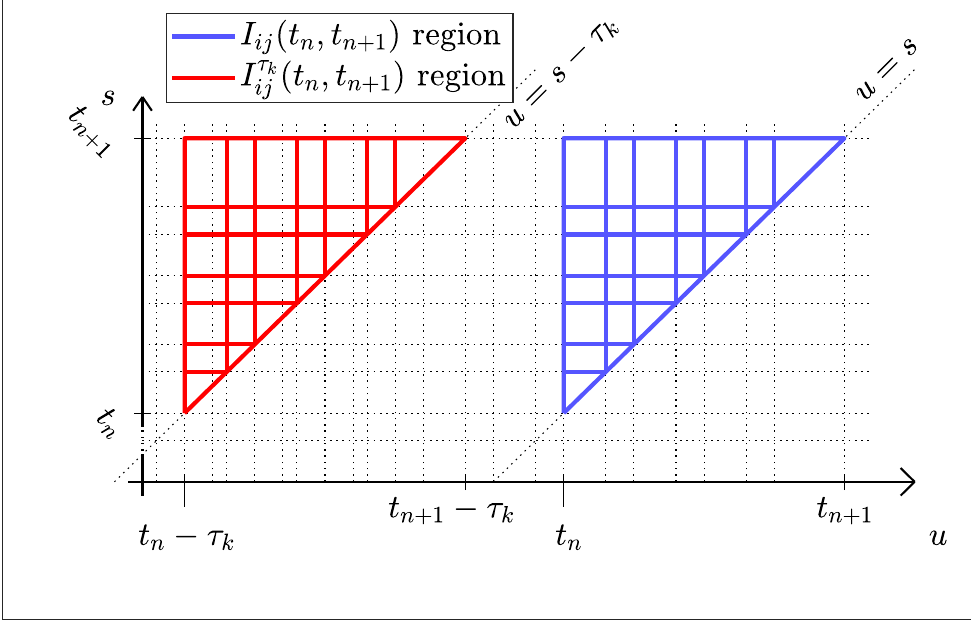}
        \caption{The ARTM includes corner points from the integration region of $I_{ij}^{\tau_k}(t_n,t_{n+1})$.}
        \label{Figure3B}
    \end{subfigure}
    \caption{Augmented mesh constructions used to simulate numerical schemes with OoC $1$, c.f.\ \Cref{Figure1B,Figure3}. \Cref{Figure1C} shows the augmented time mesh, which includes all delayed times required in the Milstein scheme \eqref{MilsteinScheme}. \Cref{Figure3B} shows the augmented refined time mesh, which includes all corner points of the delayed integration region for $I_{ij}^{\tau_k}(t_n,t_{n+1})$.}
    \label{Figures1and3Improved}
\end{figure}

In \cref{Figures1and3Improved}, we revisit the examples shown in \cref{Figure1B,Figure3}, now extended to an augmented mesh and the ARTM, respectively. We recall in \cref{Figure1B}, the Milstein scheme is not simulated at times such as $t_n-t_l$. The augmented mesh enhances the initial time mesh to include these points, shown in \cref{Figure1C}. We also enhance the mesh shown in \cref{Figure3}, so that the relevant delayed times now coincide with the ARTM, in \cref{Figure3B}.


We now extend the approximations \eqref{IijandIijtaukRefinedApproximations} for indivisible SDDEs. When simulating the delayed integrals for the Milstein implementation \eqref{MilsteinScheme}, we observe that for every time point $t_n^{(l)}\in[t_n,t_{n+1}]$ in the ARTM, the delayed time $t_n^{(l)}-\tau_k$ is also in the ARTM, due to the construction of the augmented mesh in \cref{AlgorithmAugmentedMesh}, as shown in \cref{Figure3B}. Therefore we are able to use the ARTM to simulate $I_{ij}^{\tau_k}(t_n,t_{n+1})$ in the case of indivisible delays, analogously to using a refined scheme mesh for divisible delays. That is, in the case of an indivisible SDDE, we simulate $\tilde{I}_{ij}\approx I_{ij}(t_n,t_{n+1})$ and $\tilde{I}_{ij}^{\tau_k}\approx I_{ij}^{\tau_k}(t_n,t_{n+1})$ using the trapezoidal formulae
\begin{align}
\tilde{I}_{ij}      &= \sum_{l=0}^{F_n-1}\frac{\Delta W_i(t_n^{(l)},t_n^{(l+1)})\,\Delta W_j(t_n^{(l)},t_n^{(l+1)})}{2}+\sum_{l=0}^{F_n-1}\Delta W_i(t_n^{(l)},t_n^{(l+1)})\,\Delta W_j(t_n^{(l+1)},t_n^{(F_n)}),\nonumber\\
\tilde{I}_{ij}^{\tau_k}  &= \sum_{l=0}^{F_n-1}\frac{\Delta W_i^{\tau_k}(t_n^{(l)},t_n^{(l+1)})\,\Delta W_j(t_n^{(l)},t_n^{(l+1)})}{2}+\sum_{l=0}^{F_n-1}\Delta W_i^{\tau_k}(t_n^{(l)},t_n^{(l+1)})\,\Delta W_j(t_n^{(l+1)},t_n^{(F_n)}),\label{IijandIijtaukARTMApproximations}
\end{align}
where $t_n^{(0)}=t_n$, $t_n^{(F_n)}=t_{n+1}$, and $t_n^{(0)},t_n^{(1)},\ldots,t_n^{(F_n)}$ are the ARTM points in the interval $[t_n,t_{n+1}]$, satisfying $t_n^{(0)}<t_n^{(1)}<\cdots<t_n^{(F_n)}$. The delayed values $t_n^{(0)}-\tau_k,t_n^{(1)}-\tau_k,\ldots,t_n^{(F_n)}-\tau_k\in[t_n-\tau_k,t_{n+1}-\tau_k]$ are also in the ARTM, ensuring that the delayed Wiener increments are already simulated.

\begin{remark}
As we observed above, if $t_n^{(l)}\in[t_n,t_{n+1}]$ is an ARTM value, then $t_n^{(l)}-\tau_k$ is also an ARTM value. However, the converse is not true in general, and the number of ARTM points in the (scheme) time intervals $[t_n,t_{n+1}]$ is not constant. That is, the value $F_n$ in \eqref{IijandIijtaukARTMApproximations} may vary with $n$. This is also shown in \cref{Figure3B}, where there are more ARTM points in $[t_n-\tau_k,t_{n+1}-\tau_k]$ than in $[t_n,t_{n+1}]$.
\end{remark}

\begin{theorem}\label{AugmentedMilsteinOoC1}
\textcolor{black}{Suppose $Y=(Y_n)_{n=-p,\ldots,N}$ is the Milstein approximation \eqref{MilsteinScheme} of the SDDE \eqref{Equation1}, with the integrals $I_{ij}(t_n,t_{n+1})$ and $I_{ij}^{\tau_k}(t_n,t_{n+1})$ replaced by the approximations \eqref{IijandIijtaukARTMApproximations}. That is, for $n=0,\ldots,N-1$, $Y$ satisfies
\begin{align}Y_{n+1}&=Y_n+[A_0Y_n+f(t_n,Y_n,Y_n^{\tau_1},\ldots,Y_n^{\tau_K})]\,h_n+\sum_{j=1}^m[A_jY_n+g_j(t_n,Y_n,Y_n^{\tau_1},\ldots,Y_n^{\tau_K})]\,\Delta W_j(t_n,t_{n+1})\nonumber\\
                &\;+\sum_{i=1}^m\sum_{j=1}^m\big[A_j+\nabla_xg_j(t_n,Y_n,Y_n^{\tau_1},\ldots,Y_n^{\tau_K})\big]\,\big[A_iY_n+g_i(t_n,Y_n,Y_n^{\tau_1},\ldots,Y_n^{\tau_K})\big]\,\tilde{I}_{ij}\label{AugmentedMilsteinScheme}\\
                &\;+\sum_{k=1}^K\mathbb{I}(t_n\geqslant \tau_k)\sum_{i=1}^m\sum_{j=1}^m\nabla_{x_{\tau_k}}g_j(t_n,Y_n,Y_n^{\tau_1},\ldots,Y_n^{\tau_K})\,\left[A_iY_n^{\tau_k}+g_i(t_n^{\tau_k},Y_n^{\tau_k},Y_n^{\tau_1,\tau_k},\ldots,Y_n^{\tau_K,\tau_k})\right]\,\tilde{I}_{ij}^{\tau_k}.\nonumber\end{align}
Furthermore, suppose \cref{AlgorithmAugmentedMesh} has been applied twice, to both construct the nonnegative scheme time points $(t_n)_{n=0,\ldots,N}$ using a fixed initial step size $h_\mathrm{initial}$, and also to construct an ARTM $(t_n^{(l)})_{l=0,\ldots,F_n,\,n=0,\ldots,N^\mathrm{R}}$ (where $F_1+\cdots+F_{N^\mathrm{R}}$ is the number of positive time points in the ARTM) from a refined initial step size $h^\mathrm{R}_\mathrm{initial}$ that divides $h_\mathrm{initial}$. There exists a constant $C>0$ such that
\begin{equation}\bigg(\mathbb{E}\bigg[\max_{n=-p,\ldots,N}\norm{X(t_n)-Y_n}^2\bigg]\bigg)^{1/2}\leqslant C\max_{n=-p,\ldots,N-1}\{t_{n+1}-t_n\}.\label{OoC1AgumentedMilstein}\end{equation}}
\end{theorem}

\begin{pf}
\textcolor{black}{For each $n=0,\ldots,N-1$, we choose a constant
\[C^{(n)}\geqslant\dfrac{\max_{l=0,\ldots,F_n-1}\big\{t_n^{(l+1)}-t_n^{(l)}\big\}}{4\,(t_{n+1}-t_n)^2},\]
and use \eqref{IijandIijtaukARTMApproximations} to compute
\begin{align}
\mathbb{E}&\big[\big|I_{ij}(t_n,t_{n+1})-\tilde{I}_{ij}\big|^2\big]=\mathbb{E}\bigg[\bigg|\sum_{l=0}^{F_n-1}\bigg(I_{ij}(t_n^{(l)},t_n^{(l+1)})-\frac{\Delta W_i(t_n^{(l)}\!,t_n^{(l+1)})\,\Delta W_j(t_n^{(l)}\!,t_n^{(l+1)})}{2}\bigg)\bigg|^2\bigg]\nonumber\\
&\leqslant\!\sum_{l=0}^{F_n-1}\!\int_{t_n^{(l)}}^{t_n^{(l+1)}}\!\mathbb{E}\bigg[\bigg|\int_{t_n^{(l)}}^s\,\mathrm{d}W_i(u)-\frac{\Delta W_i(t_n^{(l)}\!,t_n^{(l+1)})}{2}\bigg|^2\bigg]\,\mathrm{d}s
\leqslant\!\sum_{l=0}^{F_n-1}\!\int_{t_n^{(l)}}^{t_n^{(l+1)}}\!\bigg(\mathbb{E}\bigg[\bigg|\frac{W_i(s)-W_i(t_n^{(l)})}{2}\bigg|^2\bigg]+\mathbb{E}\bigg[\bigg|\frac{W_i(t_n^{(l+1)})-W_i(s)}{2}\bigg|^2\bigg]\bigg)\!\,\mathrm{d}s\nonumber\\
&=\!\sum_{l=0}^{F_n-1}\!\int_{t_n^{(l)}}^{t_n^{(l+1)}}\frac{1}{4}\big((s-t_n^{(l)})+(t_n^{(l+1)}-s)\big)\,\mathrm{d}s=\frac{1}{4}\sum_{l=0}^{F_n-1}(t_n^{(l+1)}-t_n^{(l)})^2\leqslant\frac{1}{4}\max_{l=0,\ldots,F_n-1}\big\{t_n^{(l+1)}-t_n^{(l)}\big\}\,\sum_{l=0}^{F_n-1}(t_n^{(l+1)}-t_n^{(l)})\nonumber\\
&=\frac{1}{4}\max_{l=0,\ldots,F_n-1}\big\{t_n^{(l+1)}-t_n^{(l)}\big\}\,(t_{n+1}-t_n)\leqslant C^{(n)}(t_{n+1}-t_n)^3,\nonumber\\
\mathbb{E}&\big[\big|I_{ij}^{\tau_k}(t_n,t_{n+1})-\tilde{I}_{ij}^{\tau_k}\big|^2\big]\leqslant C^{(n)}(t_{n+1}-t_n)^3,\label{ARTMIntergralsErrors}
\end{align}
for each $i,j=1,\ldots,m$, $i\neq j$, and $k=1,\ldots,K$. Letting $\overline{Y}=(\overline{Y}(t))_{t\in[-\tau,T]}$ be the continuous Milstein scheme from \cref{TheoremKloedenShardlow} (applied to the semilinear SDDE \eqref{Equation1}), defined on the augmented mesh $(t_n)_{n=0,\ldots,N}$ when $t\geqslant0$ (and $\overline{Y}(t)=\phi(t)$ when $t<0$), we use \eqref{ARTMIntergralsErrors} to find
\begin{align*}
\big\Vert \mathbb{E}[\overline{Y}(t_{n+1})|\overline{Y}(t_n)=x]-\mathbb{E}[Y_{n+1}|Y_n=x]\big\Vert &\leqslant C_1\,(t_{n+1}-t_n)^2,\\
\mathbb{E}\big[\norm{\overline{Y}(t_{n+1})-Y_{n+1}}^2|\overline{Y}(t_n)=Y_n=x\big]&=\mathbb{E}\bigg[\bigg|\sum_{i=1}^m\sum_{j=1}^m\big(I_{ij}(t_n,t_{n+1})-\tilde{I}_{ij}\big)+\sum_{k=1}^K\mathbb{I}(t_n\geqslant\tau_k)\sum_{i=1}^m\sum_{j=1}^m\big(I_{ij}^{\tau_k}(t_n,t_{n+1})-\tilde{I}^{\tau_k}_{ij}\big)\bigg|^2\bigg]\\
&\quad\leqslant (K+1)m^2C^{(n)}\,(t_{n+1}-t_n)^3\leqslant C_2\,(t_{n+1}-t_n)^3,
\end{align*}
for constants $C_1>0$ and $C_2=(K+1)m^2\max_{n=-p,\ldots,N-1}\{C^{(n)}\}$, so that
\begin{equation}\mathbb{E}\bigg[\max_{n=-p,\ldots,N}\norm{\overline{Y}(t_n)-Y_n}^2\bigg]\leqslant C_3\max_{n=-p,\ldots,N-1}\{t_{n+1}-t_n\}^2,\label{MSOoCofAugmentedMilstein}\end{equation}
for constant $C_3=\max\{C_1,C_2\}$, by \cref{GeneralTheoremofMilsteinforSDDEs}. Furthermore, \cref{TheoremKloedenShardlow} shows
\begin{equation}\mathbb{E}\bigg[\max_{n=-p,\ldots,N}\norm{X(t_n)-\overline{Y}(t_n)}^2\bigg]\leqslant C_4\max_{n=-p,\ldots,N-1}\{t_{n+1}-t_n\}^2,\label{C4}\end{equation}
for constant $C_4>0$. Combining \eqref{MSOoCofAugmentedMilstein} with \eqref{C4} and letting $C=C_3+C_4$ gives
\begin{align*}
\mathbb{E}\bigg[\max_{n=-p,\ldots,N}\norm{X(t_n)-Y_n}^2\bigg]&=\mathbb{E}\bigg[\max_{n=-p,\ldots,N}\norm{X(t_n)-\overline{Y}(t_n)}^2\bigg]+\mathbb{E}\bigg[\max_{n=-p,\ldots,N}\norm{\overline{Y}(t_n)-Y_n}^2\bigg]\leqslant C\max_{n=-p,\ldots,N-1}\{t_{n+1}-t_n\}^2,
\end{align*}
from which \eqref{OoC1AgumentedMilstein} holds.
\hfill$\blacksquare$}
\end{pf}

\begin{remark}
\textcolor{black}{In our proof of \cref{AugmentedMilsteinOoC1}, we have appealed to the result \cref{TheoremKloedenShardlow}, assuming global Lipschitz conditions due to continuous differentiability (see \cite[page 183]{KloedenShardlow2012}). However, while producing \cref{LargerTExample}, we also found that replacing the parameters \eqref{Example1Parameters} with functions satisfying only local and not global Lipschitz conditions resulted in Milstein schemes (on the augmented mesh) that were numerically unstable. Future work may address this, such as by introducing a tamed Milstein scheme defined on the augmented mesh.}
\end{remark}

\begin{remark}
\textcolor{black}{Alternative approaches for simulating iterated stochastic integrals include Fourier–series–based methods, such as the expansion introduced by Wiktorsson \cite{Wiktorsson2001}. Similarly to the case of non–delayed SDEs, these methods are typically built upon
trapezoidal–type decompositions of the integration region, supplemented by additional Fourier–series expansions and auxiliary random variables. These methods may be extended to SDDEs, but the objective of this work is to introduce the method of the augmented mesh as a framework for simulating SDDE solutions when the delays are indivisible. We therefore restrict attention to direct trapezoidal approximations, without the additional algorithmic complexity associated with
Fourier–based implementations.}
\end{remark}

\section{Examples using the augmented mesh}\label{Section7}

In this section, we apply the augmented time mesh to simulate numerical solutions of indivisible SDDEs, \textcolor{black}{for three examples}. In the first example, we reuse the equation and parameters from \cref{DivisibleExample} but vary the delay values, using the same delays as in \cref{LIExample}. \textcolor{black}{We also present a selection of the computation times for the first example, along with discussion of when the method of the augmented mesh becomes practically advantageous over other methods.} In the second example, we fix all parameters and delays, to investigate how the accuracy of the numerical solution evolves as the simulation time increases. For both of the first two examples, we demonstrate scheme performance using error graphs adapted to the indivisible setting. \textcolor{black}{We then give a third practical example, showing the applicability of augmented meshes to realistic models, such as that described in \cref{FinancialExample}.}

In the divisible setting, error graphs such as those in \cref{ErrorGraphs1} show numerical errors against constant step sizes $h$, enabling simple interpretation such as halving the step sizes (selecting $h=2^{-3},2^{-4},\ldots$) corresponding to halving the errors for schemes with OoC $1$. However, for indivisible delays, the augmented mesh introduces varying step sizes $h_n$, requiring an alternative visualisation. Rather than plotting the errors against the step sizes, we instead plot the errors against the initial step sizes. Additionally, the linear reference curves (in the $\log_2$-$\log_{10}$ scale) are only meaningful when step sizes divide (halve in the case of \cref{ErrorGraphs1}) uniformly, so we omit these when the delays are indivisible.

\begin{example}\label{IndivisibleExample}
This example replicates \cref{DivisibleExample}, but varies the delay values $\tau_1$ and $\tau_2$ to test the robustness of the numerical schemes under indivisible conditions. We again simulate numerical solutions of equation \eqref{Equation1} with $m=2$, $d=2$, $T=4$, $\phi$ as given by \eqref{phit}, while $A_0,A_1,A_2,f,g_1,g_2$ are given by \eqref{Example1Parameters}. We select the delay pairs $(\tau_1,\tau_2)=(1,1/2),\,(1,\pi/4),\,(\exp(2)/5,\pi/4)$, and show the error graphs of the numerical schemes at time $t_n=T$, simulated for initial step sizes 
 $h_\mathrm{initial}=2^{-1}\ldots,2^{-10}$. 
 These simulations are produced using initial refined step size $h^\mathrm{R}_\mathrm{initial}=2^{-13}$, while the error graphs are simulated with $n_\mathrm{t}=1,000$ trials.

\begin{figure}[ht]
    \centering
    \begin{subfigure}[t]{0.32\linewidth}
        \centering
        \includegraphics[width=\linewidth,trim=0 0 0 0,clip]{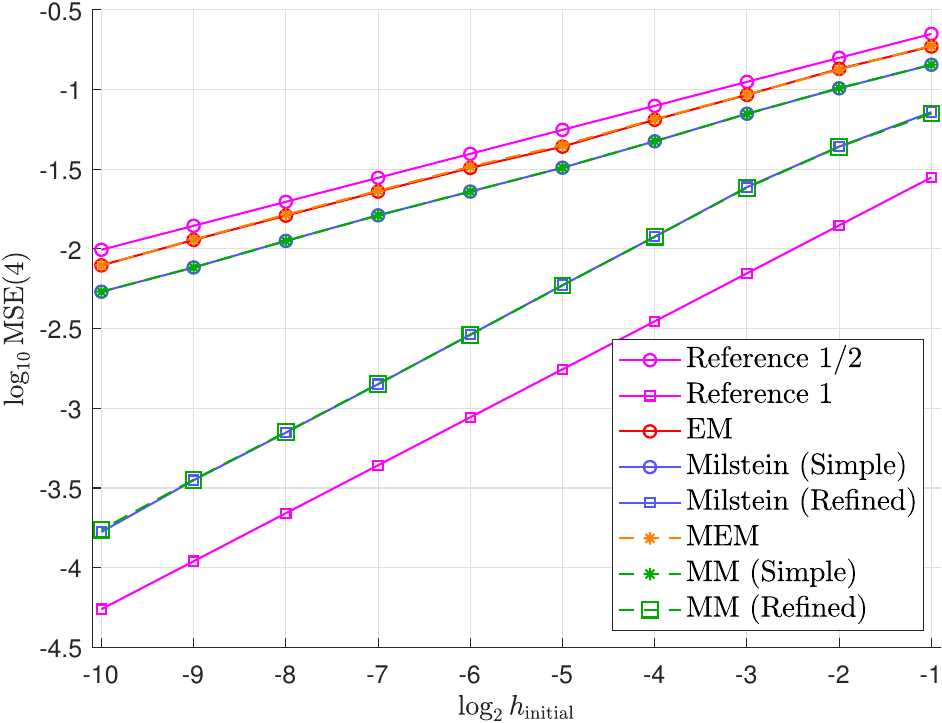}
        \caption{Error graphs with delays $\tau_1=1$ and $\tau_2=1/2$.}
        \label{Example2A}
    \end{subfigure}
    \hfill
    \begin{subfigure}[t]{0.32\linewidth}
        \centering
        \includegraphics[width=\linewidth,trim=0 0 0 0,clip]{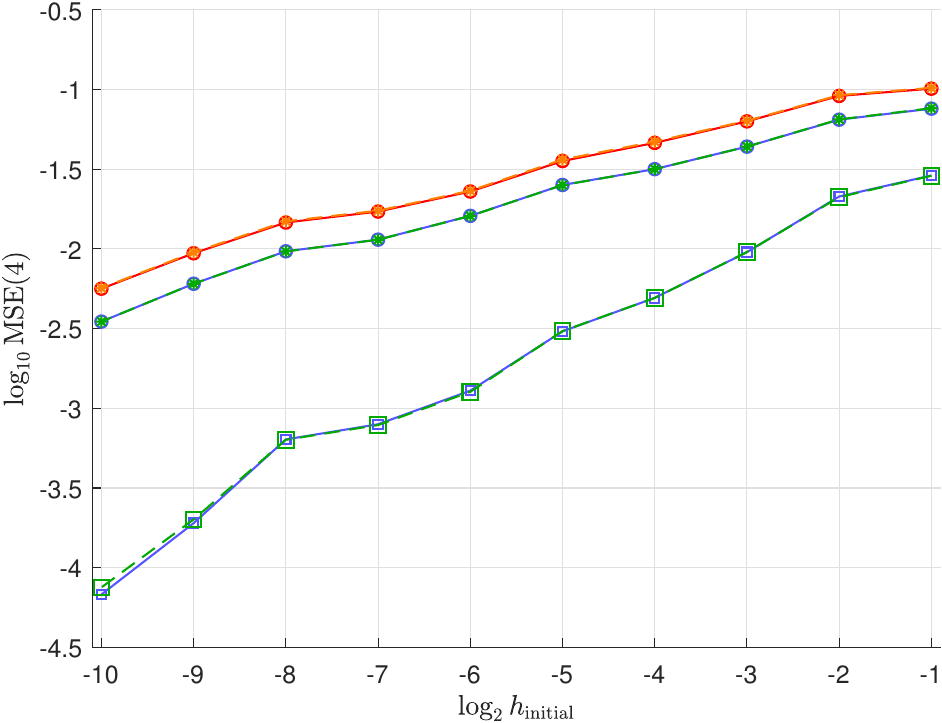}
        \caption{Error graphs with delays $\tau_1=1$ and $\tau_2=\pi/4$.}
        \label{Example2B}
    \end{subfigure}
    \hfill
    \begin{subfigure}[t]{0.32\linewidth}
        \centering
        \includegraphics[width=\linewidth,trim=0 0 0 0,clip]{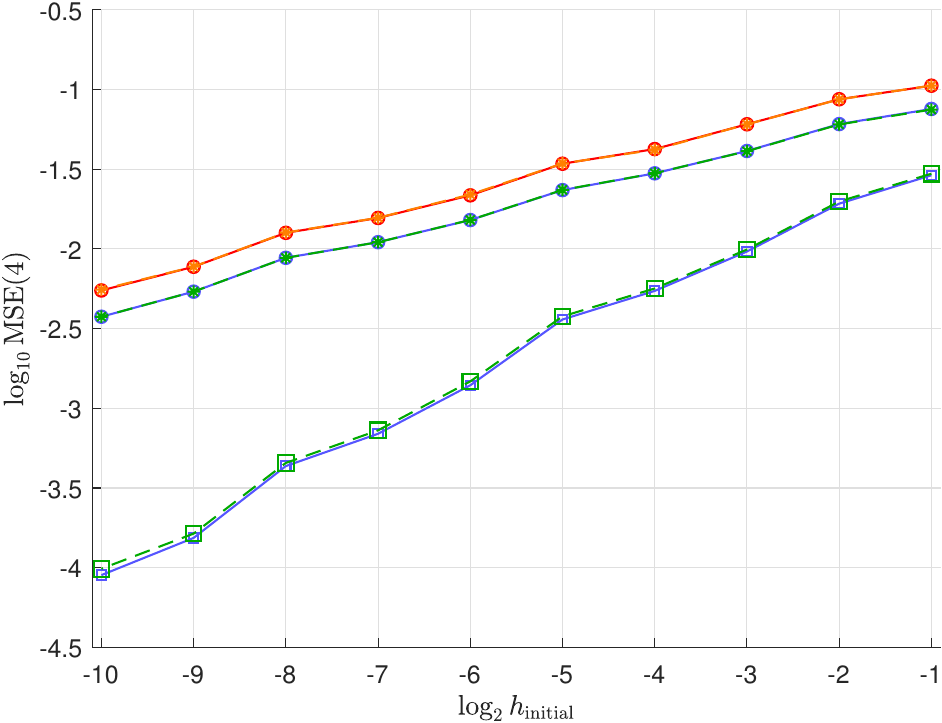}
        \caption{Error graphs with $\tau_1=\exp(2)/5$ and $\tau_2=\pi/4$.}
        \label{Example2C}
    \end{subfigure}
    \caption{Mean-square error graphs for the numerical solutions of equation \eqref{Equation1} applied on augmented meshes, with the parameters given in \cref{DivisibleExample}, using three different pairs of delays $\tau_1$ and $\tau_2$. Each scheme is simulated on an augmented mesh constructed from an initial step size $h_\mathrm{initial}$.}
    \label{ErrorGraphs2}
\end{figure}

\cref{ErrorGraphs2} shows the error graphs for the numerical schemes simulated on augmented time meshes. For given initial step size $h_\mathrm{initial}$, each scheme is simulated on the augmented mesh constructed from $h_\mathrm{initial}$, with error $\mathrm{MSE}(4)$ calculated at time $4$. In the case of \cref{Example2A}, the delays are divisible, so the augmented mesh coincides with the initial time mesh, and in this case, the error graphs are linear on the $\log_2$-$\log_{10}$ scale. Furthermore, \cref{Example1A} is identical to \cref{Example2A} because the delays are identical in these cases. \cref{Example2B} shows the case where only one delay is indivisible (with respect to each initial step size), while in \cref{Example2C}, both delays are indivisible (with every initial step size). In both of \cref{Example2B,Example2C}, the indivisible delays cause the step sizes of the schemes to be no longer constant, corresponding to nonlinear error graphs. In either of these latter cases, however, the schemes with OoC $1$ have errors significantly less than the errors of the schemes with OoC $1/2$.

\end{example}

\textcolor{black}{\Cref{Table3} shows total computation times for the EM LI scheme applied in \cref{LIExample}, together with the times for the numerical schemes from \cref{IndivisibleExample}. To produce these times, all simulations were performed using MATLAB~R2023b on a computer with a 13th Gen Intel Core i5-13600K processor and 16\,GB RAM. Each entry shows the total computation time across $n_\mathrm{t}=1,000$ trials, for each scheme. Times are shown for both divisible delays $(\tau_1,\,\tau_2)=(1,\,1/2)$ and indivisible delays $(\tau_1,\,\tau_2)=(\exp(2)/5,\,\pi/4)$. In the case of divisible delays, all schemes are simulated with a fixed step size $h$. For indivisible delays, the EM LI scheme is simulated using $h_\mathrm{initial}$ as a fixed step size, while all other schemes are simulated on the augmented time mesh constructed from initial step size $h_\mathrm{initial}$. To isolate the main contributors to computational expense, \Cref{Table3} lists times for key components separately from the scheme evaluations. Component times include the refined integrals ($\tilde I_{0j}$, $\tilde I_{12}$, $\tilde I_{ij}^{\tau_k}$) computed on the ARTM using \eqref{IijandIijtaukARTMApproximations}, as well as the Magnus exponents ($M_n^{[1]}$ and $M_n^{[2]}$, using both simple and refined integral approximations). Scheme times include only the evaluations of \eqref{EMScheme}, \eqref{MilsteinScheme}, \eqref{MEMScheme}, \eqref{MMScheme}, excluding component costs.}

\begin{table}
\centering
  \begin{tabular}{|c|c|c|c|c|c|c|c|c|c|c|}
    \hline
    \multicolumn{11}{|c|}{\vphantom{$\displaystyle{\int}$}\textbf{Comparison of Computation Times for SDDE Schemes, Divisible vs Indivisible}}\\
    \hline
    \multirow{2}{*}{Scheme/Component}&
      \multicolumn{5}{|c|}{Step Size $h$ (Divisible Case)}&\multicolumn{5}{|c|}{Initial Step Size $h_\mathrm{initial}$ (Indivisible Case)} \\
&$2^{-2}$&$2^{-3}$&$2^{-4}$&$2^{-5}$&$2^{-6}$&$2^{-2}$&$2^{-3}$&$2^{-4}$&$2^{-5}$&$2^{-6}$ \\
    \hline
    EM LI&0.021&0.036&0.049&0.094&0.185&0.028&0.029&0.065&0.109&\cellcolor[RGB]{150,255,255}0.177\\
    \hline
    EM&0.009&0.016&0.035&0.063&0.141&0.031&0.048&0.085&0.169&0.324\\
    \hline
    Milstein (Simple)&0.014&0.029&0.059&0.123&0.240&0.051&0.086&0.159&0.302&0.583\\
    \hline
    Milstein (Refined)&0.015&0.029&0.059&0.119&0.241&\cellcolor[RGB]{255,100,100}0.050&0.085&0.153&0.286&0.564\\
    \hline
    MEM&0.009&0.017&0.033&0.066&0.145&0.027&0.045&0.084&0.155&0.311\\
    \hline
    MM (Simple)&0.020&0.040&0.079&0.161&0.324&0.071&0.119&0.219&0.416&0.807\\
    \hline
    MM (Refined)&0.019&0.039&0.077&0.155&0.319&0.072&0.121&0.216&0.407&0.786\\
    \hline
    $\vphantom{\sum^\int_l}\tilde I_{0j}(t_n,t_{n+1})$, $j=1,2$&0.016&0.017&0.018&0.021&0.024&\cellcolor[RGB]{255,100,100}0.042&0.043&0.047&0.051&0.061\\
    \hline
    $\vphantom{\sum^\int_l}\tilde I_{12}(t_n,t_{n+1})$&0.009&0.010&0.011&0.013&0.018&\cellcolor[RGB]{255,100,100}0.024&0.025&0.029&0.035&0.046\\
    \hline
    $\vphantom{\sum^\int_l}\tilde I_{ij}^{\tau_k}(t_n,t_{n+1})$, $i,j,k=1,2$&0.032&0.032&0.035&0.039&0.048&\cellcolor[RGB]{255,100,100}0.047&0.049&0.053&0.061&0.077\\
    \hline
    $\vphantom{\sum^\int_l}M_n^{[1]}(t_n,t_{n+1})$&0.029&0.056&0.112&0.223&0.443&0.114&0.188&0.338&0.629&1.215\\
    \hline
    $\vphantom{\sum^\int_l}M_n^{[2]}(t_n,t_{n+1})$ (Simple)&0.024&0.046&0.091&0.183&0.365&0.090&0.149&0.268&0.499&0.963\\
    \hline
    $\vphantom{\sum^\int_l}M_n^{[2]}(t_n,t_{n+1})$ (Refined)&0.027&0.053&0.104&0.208&0.414&0.101&0.168&0.300&0.558&1.076\\
    \hline
    Reference Solution&\multicolumn{5}{|c|}{35.297}& \multicolumn{5}{|c|}{92.787} \\
    \hline
    Simulation Time&\multicolumn{5}{|c|}{66.844}& \multicolumn{5}{|c|}{206.549} \\
    \hline
  \end{tabular}
  \caption{\textcolor{black}{Computation times (in minutes) for $n_\mathrm{t}=1,000$ simulations of our numerical schemes applied to a $2$-dimensional SDDE until the terminal time $T=4$, comparing times for divisible delays $(\tau_1,\,\tau_2)=(1,\,1/2)$, and indivisible delays $(\tau_1,\,\tau_2)=(\exp(2)/5,\,\pi/4)$. Times shown in the columns beneath $h$ correspond to divisible delays, while times beneath $h_{\mathrm{initial}}$ correspond to indivisible delays. In the case of indivisible delays, the EM LI scheme uses $h_\mathrm{initial}$ as a fixed step size, while all other schemes are constructed on an augmented time mesh with $h_\mathrm{initial}$ as the initial step size. Times listed alongside the numerical schemes (EM LI, EM, Milstein (Simple), Milstein (Refined), MEM, MM (Simple), and MM (Refined)) state the computation times of those schemes without including the times to evaluate particular key components ($\tilde I_{0j}(t_n,t_{n+1})$, $\tilde I_{12}(t_n,t_{n+1})$, $\tilde I_{ij}^{\tau_k}(t_n,t_{n+1})$, for $i,j,k=1,2$, $M_n^{[1]}(t_n,t_{n+1})$, $M_n^{[2]}(t_n,t_{n+1})$ (Simple), and $M_n^{[2]}(t_n,t_{n+1})$ (Refined)), while the times to evaluate key components are listed separately. The computation time to simulate the reference solutions is also listed, as is the total simulation time, for both the divisible and indivisible cases.}}\label{Table3}
\end{table}

\textcolor{black}{\cref{Table3} confirms that the augmented-mesh approach incurs increased computational cost in exchange for higher accuracy. As indicated by bounds such as \eqref{UpperBoundonTc}, the size of the augmented mesh grows quadratically with the number of delay multiples, leading to a corresponding increase in the number of scheme evaluations and integral computations. \cref{Table3} also illustrates that when delays are relatively large compared with the terminal time, the method of the augmented mesh can be computationally advantageous. For example, the EM LI scheme with step size $2^{-6}$ attains an error of approximately $10^{-1.6}$ (see \cref{Example1LIC}), with total computation time 0.177 minutes (highlighted in cyan, in \cref{Table3}). Comparable accuracy is achieved by the refined Milstein scheme on the augmented mesh (see \cref{Example2C}), with a total computation time of 0.163 minutes (highlighted in red). In this setting, the augmented-mesh Milstein scheme is more efficient than the EM LI scheme is, with similar levels of accuracy.}

We now give a second example, demonstrating the effect of including additional mesh times from a larger terminal time $T$. This investigation is motivated by the observation that in \cref{IndivisibleExample}, while simulating the numerical solutions until time $t_n=T=4$, the density of points in the augmented mesh decreases in the later time windows such as $[3.5,4]$. For example, when the initial step size is $h_\mathrm{initial}=1/2$ (with delays satisfying \eqref{Maxh} such as $(\tau_1,\tau_2)=(\exp(2)/5,\pi/4)\approx(1.48,0.79)$), then the only augmented mesh points between $3.5$ and $4$ are those points themselves, while if $T=5$ instead, then the augmented mesh includes the additional point $T-\tau_1\approx 3.52$.

\begin{example}\label{LargerTExample}
We now replicate \cref{DivisibleExample} with parameters \eqref{Example1Parameters}, initial process \eqref{phit}, and again with $d=m=2$, but now until time $T=6$. We use the delay values $(\tau_1,\tau_2)=(1,\pi/4)$, and observe the errors of the numerical schemes at initial step sizes $h_\mathrm{initial}=2^{-1},\ldots,2^{-10}$, using initial refined step size $h_\mathrm{initial}^\mathrm{R}=2^{-13}$ to simulate the reference solutions. We show the error graphs at six different observation times, simulated with $n_\mathrm{t}=2,500$ trials.
\begin{enumerate}
  \item Time $t_n=1/2$ shows the errors after one step of the largest $h_\mathrm{initial}$ value.
  \item Time $t_n=\tau_1\approx1.48$ shows the errors after the second Bellman interval.
  \item Time $t_n=2\tau_2\approx1.57$ shows the errors at another chosen time not on the initial time mesh.
  \item Time $t_n=2$ shows the errors at time $T/3$.
  \item Time $t_n=4$ shows the errors at time $2T/3$.
  \item Time $t_n=6$ shows the errors at time $T$.
\end{enumerate}
The first three of these times are chosen to demonstrate that although the error graphs are nonlinear, they show distinct differences between the OoC-$1/2$ and OoC-$1$ schemes, regardless of the time $t_n$ at which the errors are observed. The final three of these times are chosen to show the steadying behaviour as time increases towards $T$.

\begin{figure}[ht]
    \centering
    \begin{subfigure}[t]{0.32\linewidth}
        \centering
        \includegraphics[width=\linewidth,trim=0 0 0 0,clip]{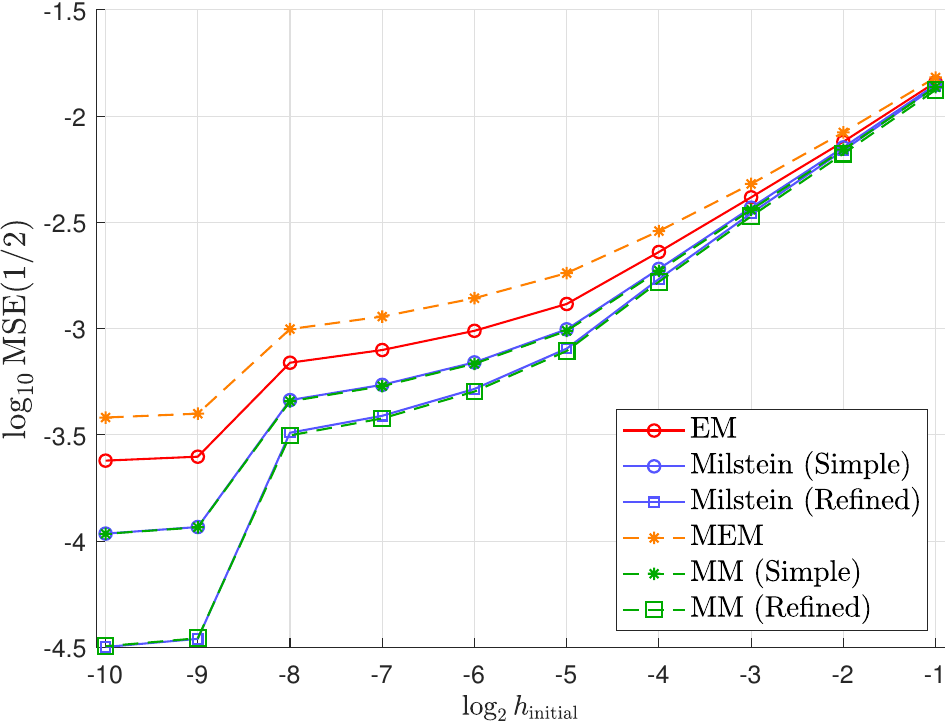}
        \caption{Error graphs at time $t_n=1/2$.}
        \label{Example3A}
    \end{subfigure}
    \hfill
    \begin{subfigure}[t]{0.32\linewidth}
        \centering
        \includegraphics[width=\linewidth,trim=0 0 0 0,clip]{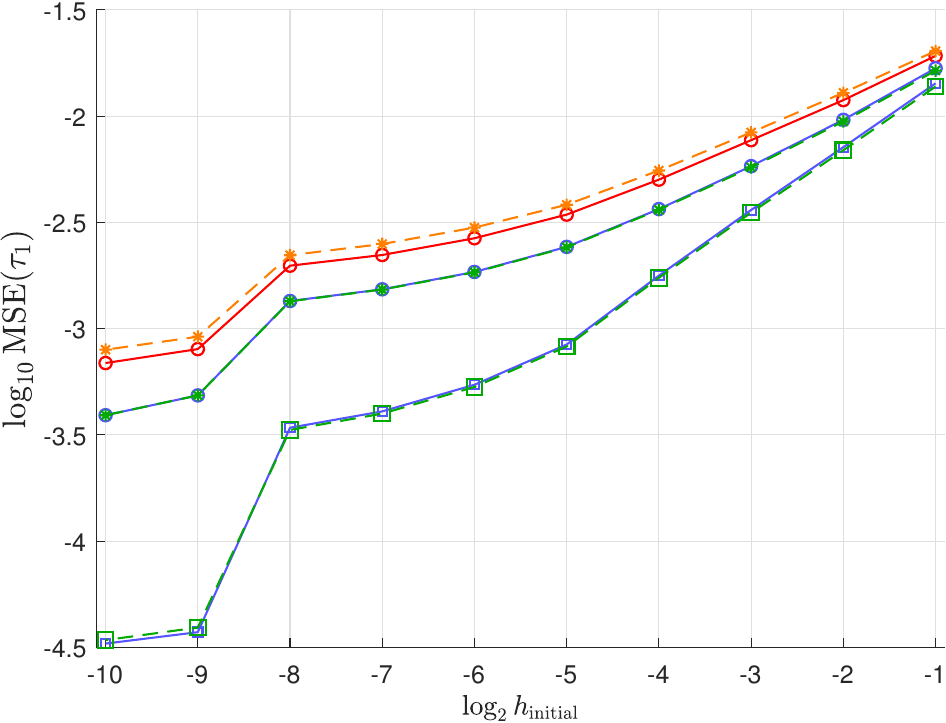}
        \caption{Error graphs at time $t_n=\tau_1\approx1.48$.}
        \label{Example3B}
    \end{subfigure}
    \hfill
    \begin{subfigure}[t]{0.32\linewidth}
        \centering
        \includegraphics[width=\linewidth,trim=0 0 0 0,clip]{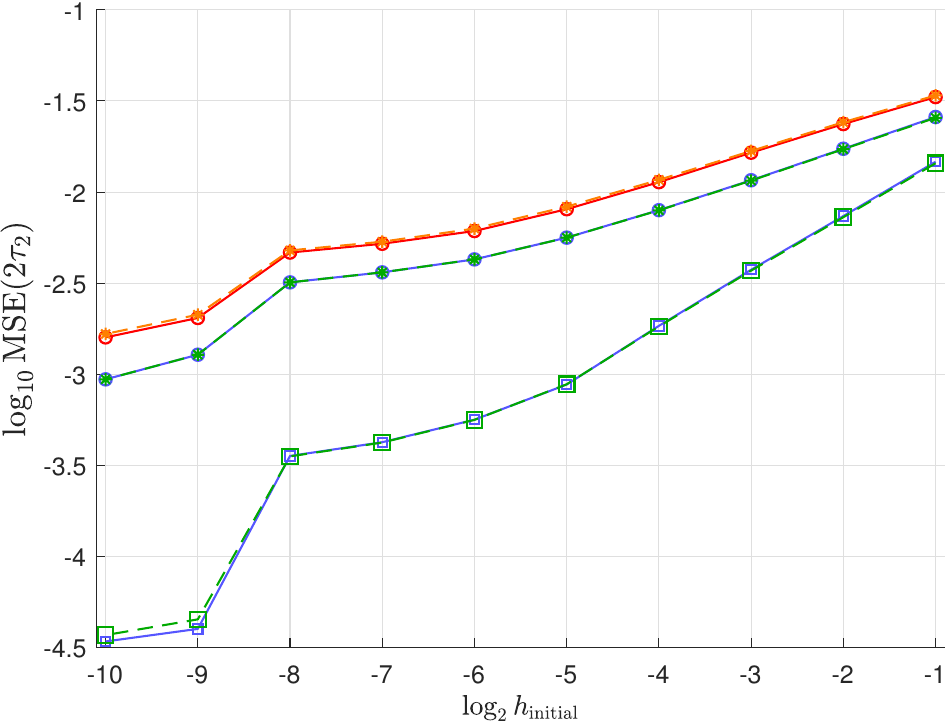}
        \caption{Error graphs at time $t_n=2\tau_2\approx1.57$.}
        \label{Example3C}
    \end{subfigure}
    \par\bigskip
    \begin{subfigure}[t]{0.32\linewidth}
        \centering
        \includegraphics[width=\linewidth,trim=0 0 0 0,clip]{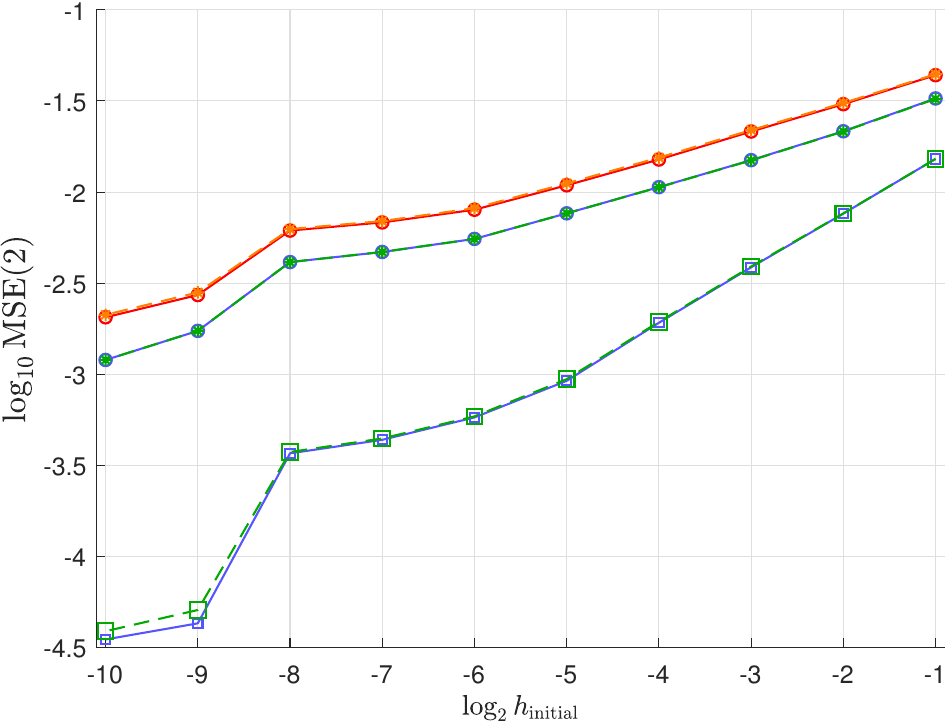}
        \caption{Error graphs at time $t_n=T/3=2$.}
        \label{Example3D}
    \end{subfigure}
    \hfill
    \begin{subfigure}[t]{0.32\linewidth}
        \centering
        \includegraphics[width=\linewidth,trim=0 0 0 0,clip]{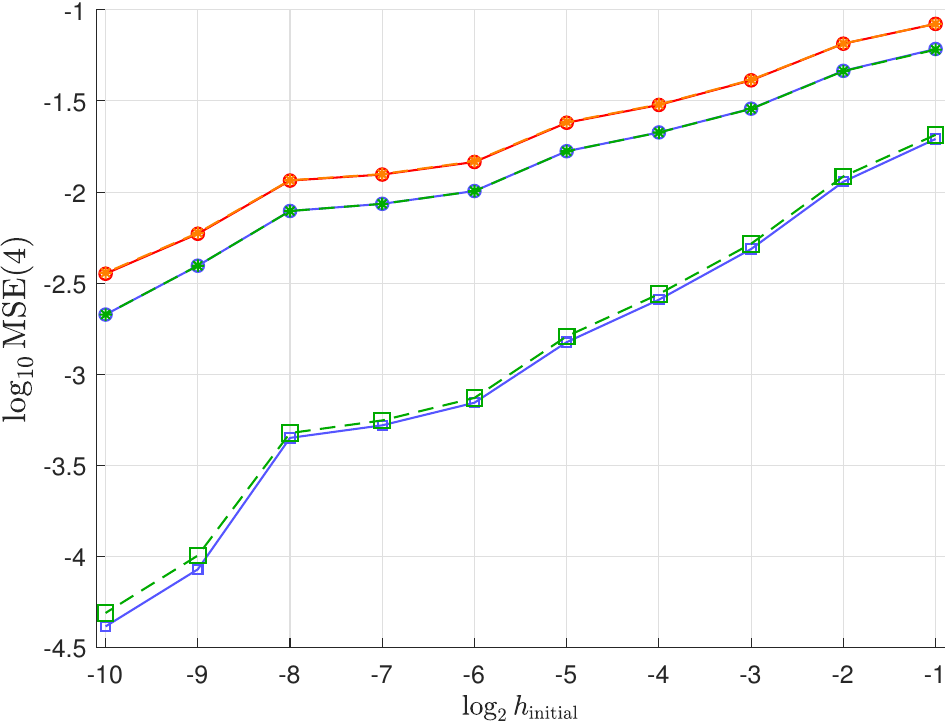}
        \caption{Error graphs at time $t_n=2T/3=4$.}
        \label{Example3E}
    \end{subfigure}
    \hfill
    \begin{subfigure}[t]{0.32\linewidth}
        \centering
        \includegraphics[width=\linewidth,trim=0 0 0 0,clip]{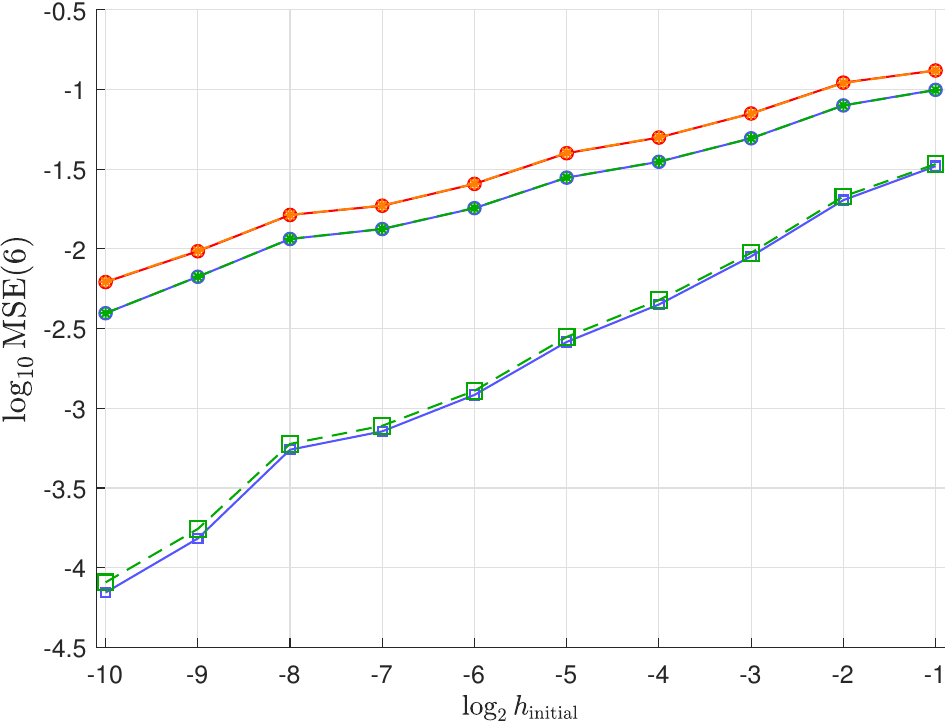}
        \caption{Error graphs at time $t_n=T=6$.}
        \label{Example3F}
    \end{subfigure}
    \caption{Mean-square error graphs for the numerical solutions of equation \eqref{Equation1} using the parameters given in \cref{DivisibleExample}, with delays $\tau_1=\exp(2)/5$ and $\tau_2=\pi/4$. In the first row of these figures, the errors show distinct difference between the OoC-$1/2$ and the OoC-$1$ schemes, regardless of which time the schemes are observed at. In the second row, the simulation is continued, increasing towards the terminal time $T=6$, with the error graphs becoming increasingly close to linear.}
    \label{ErrorGraphs3}
\end{figure}

\cref{ErrorGraphs3} shows that the schemes with theoretical OoC $1$ consistently outperform those with OoC $1/2$ across all observation times. The trend of improved accuracy (achieved by using the (refined) Milstein and MM schemes) remains robust across observation times such as $t_n=1/2,\,\tau_1,\,2\tau_2$. As the simulation progresses toward the terminal time $T=6$, the errors remain stable and the error graphs become increasingly linear. Notably, both \cref{Example3E} and \cref{Example2B} display errors with the same delays, at time $t_n=4$, but the former uses a greater terminal time ($T=6$) compared to the latter ($T=4$). 
 Furthermore, the magnitude of the errors is less in \cref{Example3E} than in \cref{Example2B} (for example, the Milstein scheme has minimal error close to $10^{-4.5}$ in \cref{Example3E}, but closer to $10^{-4}$ in \cref{Example2B}), since increasing $T=4$ to $T=6$ produces more scheme points in the augmented mesh, thus reducing the scheme error.

\end{example}


The examples in this section thus far confirm that the augmented mesh enables the (refined) Milstein and MM schemes to attain their theoretical convergence orders, even in the presence of indivisible delays. We have shown that the method of the augmented mesh preserves the convergence orders of the schemes, across a range of delay values. We have also shown that including more points in the augmented mesh, such as is caused by including more observation times, leads to decreased errors, validating the theoretical error bound given in \cref{TheoremKloedenShardlow}. \textcolor{black}{We now give a third, practical example, utilising the method of the augmented mesh.}

\begin{example}\label{SimulatedFinancialExample}

\textcolor{black}{Inspired by \cref{FinancialExample}, we estimate the value of a financial (digital call) option $V(T)=\mathbb{E}[\mathbb{I}(X(T)>K)]=\mathbb{P}(X(T)>\kappa)$, for constant (strike price) $\kappa=1.54$, assuming $X$ is the delayed local-volatility process satisfying
\begin{align}
\mathrm{d}X(t)&=\mu X(t)\,\mathrm{d}t+\frac{3\sigma}{10}X(t-1)\,\mathrm{d}W_1(t)+\frac{7\sigma}{10}X(t-\pi/4)\,\mathrm{d}W_2(t),\quad t\in[0,T],\label{DLVM}\\
X(t)&=1,\quad t\in[-1,0].\nonumber
\end{align}
We estimate $V(T)$ by simulating the process \eqref{DLVM} $n_\mathrm{t}$ times, with the EM LI scheme, Milstein LI scheme, and the Milstein scheme on the augmented mesh (using refined integral approximations). For a selection of (fixed and initial) step sizes, we show sample averaged values for each scheme, as $n_\mathrm{t}$ increased, corresponding to estimates for $V(T)$.}

\begin{figure}[ht]
\centering
\begin{subfigure}[b]{0.23\linewidth}
  \includegraphics[width=\linewidth, trim=0 0 0 0, clip]{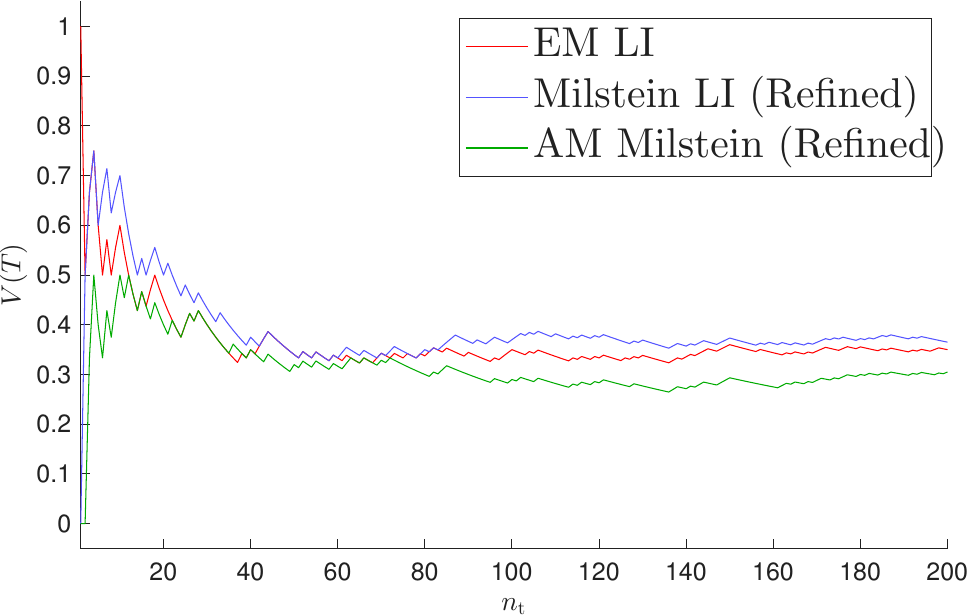}
  \caption{Step sizes $h=h_\mathrm{initial}=2^0$.}
  \label{FigureFinance1}
\end{subfigure}
\hfill
\begin{subfigure}[b]{0.23\linewidth}
  \includegraphics[width=\linewidth, trim=0 0 0 0, clip]{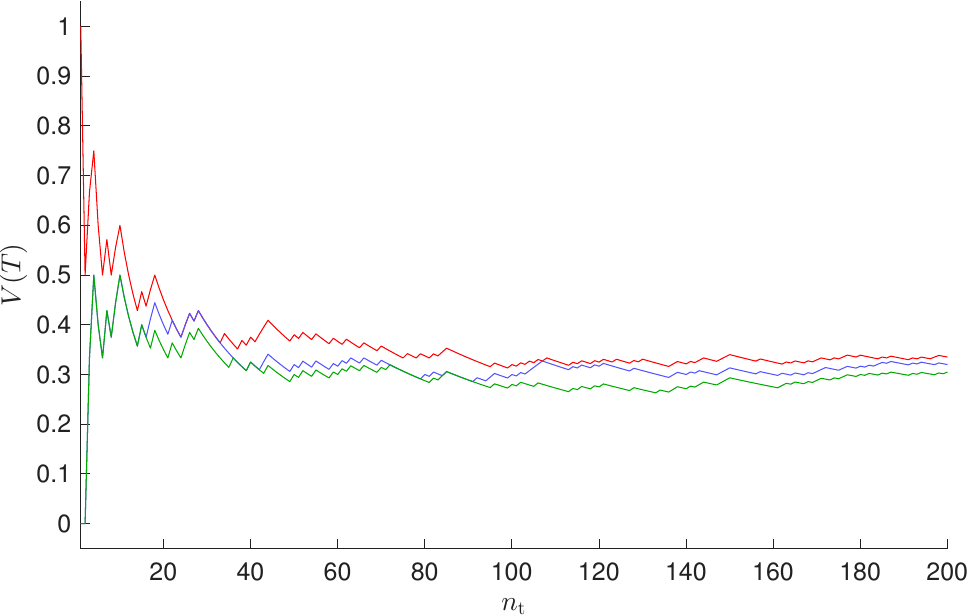}
  \caption{Step sizes $h=h_\mathrm{initial}=2^{-1}$.}
  \label{FigureFinance2}
\end{subfigure}
\hfill
\begin{subfigure}[b]{0.23\linewidth}
  \includegraphics[width=\linewidth, trim=0 0 0 0, clip]{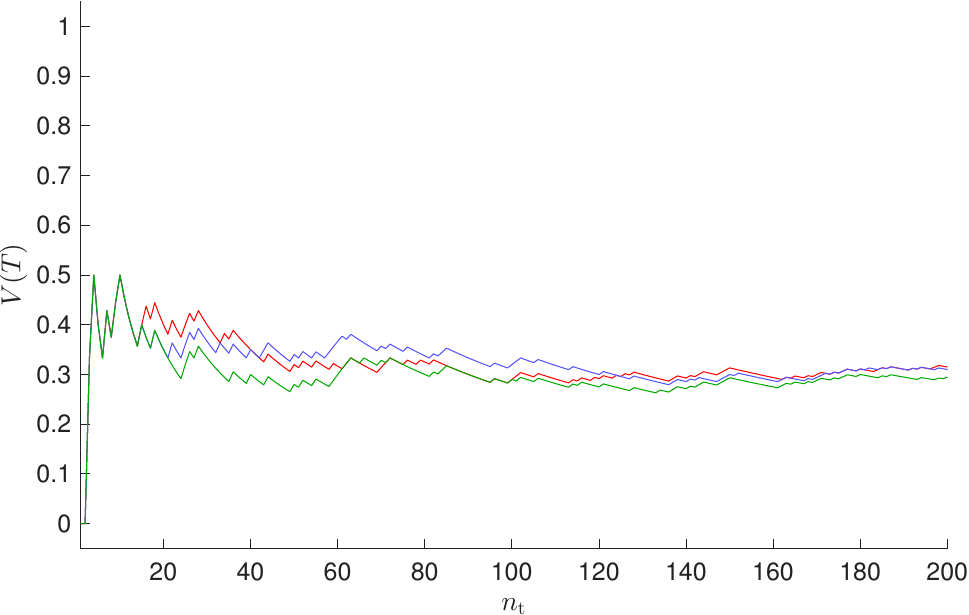}
  \caption{Step sizes $h=h_\mathrm{initial}=2^{-2}$.}
  \label{FigureFinance3}
\end{subfigure}
\hfill
\begin{subfigure}[b]{0.23\linewidth}
  \includegraphics[width=\linewidth, trim=0 0 0 0, clip]{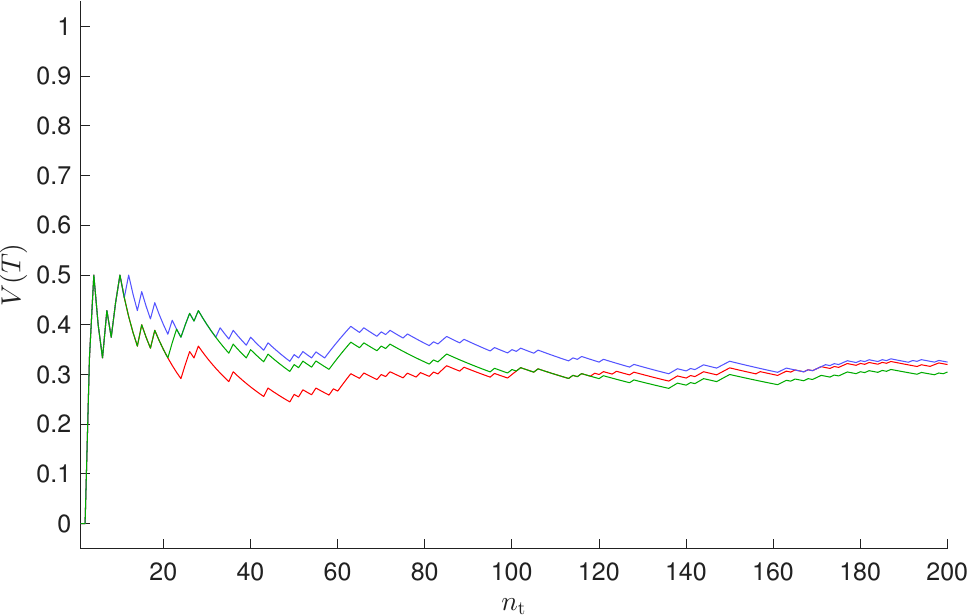}
  \caption{Step sizes $h=h_\mathrm{initial}=2^{-3}$.}
  \label{FigureFinance4}
\end{subfigure}
\caption{Estimation of the digital call option value $V(T)=\mathbb E[\mathbb I(X(T)>K)]$ for the SDDE \eqref{DLVM}. Sample–averaged estimates are shown as functions of the number of trials $n_{\mathrm t}$ for the EM LI scheme, the Milstein LI scheme, and the augmented–mesh Milstein scheme using refined integral approximations. Each subplot corresponds to a different (fixed or initial) step size $h=h_{\mathrm{initial}}$, demonstrating how the augmented mesh enables accurate option valuation for indivisible delays at coarser time resolutions.}
\label{FigureFinanceOptions}
\end{figure}

\end{example}

\section{Conclusions and discussions}\label{Section8}


This paper was motivated by the absence of numerical methods for simulating SDDE solutions with strong order of convergence $1$ in the presence of general, indivisible delays. While the work of Kloeden and Shardlow \cite{KloedenShardlow2012} (see \cref{TheoremKloedenShardlow}) guarantees the convergence order of the Milstein scheme in theory, there has been no standard implementation strategy to realize this rate in practice. To address this, we introduced the augmented time mesh as a systematic approach that enables the Milstein scheme to remain applicable and convergent even when the delays do not align with uniform step sizes.


We described the complications of implementing an OoC-$1$ scheme and then presented our solution. We first observed from the Milstein scheme \eqref{MilsteinScheme} that implementing such a scheme requires prior evaluation of scheme values at time points $t_n$, $t_n-\tau_k$, $t_n-\tau_l-\tau_k$, and so on, which necessitates a varying scheme step size when the delays are indivisible. We then extended the trapezoidal method to simulate the integrals $I_{ij}(t_n,t_{n+1})$ and $I_{ij}^{\tau_k}(t_n,t_{n+1})$ from the divisible-delay to the indivisible-delay setting. Finally, we provided numerical simulations demonstrating that the method of the augmented mesh is robust and clearly distinguishes the performance of order-$1/2$ and order-$1$ schemes.


The augmented mesh is most effective when the delays are large compared with the terminal time $T$, but as shown in \cref{Table2}, its size can grow prohibitively large in the presence of small delays. In such cases, alternative methods are possible, such as simulating delayed scheme values using linear interpolation (see \cref{Section4}). These approaches offer greater efficiency at the cost of a reduced convergence order ($1/2$ instead of $1$). The fewer delay multiples $i_k\tau_k$ that lie within $[0,T]$, the more practical the augmented mesh approach becomes.


Aside from providing a framework to simulate numerical schemes that realize convergence order $1$, our construction of the augmented mesh for SDDEs also enables the simulation of higher-order numerical solutions, as we noted in \cref{RemarkonGeneralisabilityforHigherOrderSchemes}. Future research may explore ways to improve computational efficiency, such as adaptive refinement of the augmented mesh tailored to specific delay structures. The augmented mesh is also well suited to parallel implementation, making it a promising foundation for large-scale simulations and high-dimensional SDDEs with indivisible delays. Additionally, since the augmented mesh refines any fixed-step-size mesh such as $(nh_\mathrm{initial})_{n=0,\ldots,T/h_\mathrm{initial}}$, it naturally integrates with multilevel Monte Carlo and other variance-reduction techniques. These connections offer further opportunities to advance the efficient simulation of complex stochastic delay models across scientific and engineering applications.

\section*{Acknowledgments}
This research was conducted during Mitchell Griggs' PhD studies at Queensland University of Technology (QUT), supported by the QUT Postgraduate Research Award (QUTPRA) and internal research funding provided by Kevin Burrage through the School of Mathematical Sciences.

The authors are grateful to the anonymous reviewers for valuable feedback that improved the quality of this paper.

\section*{Competing Interests}
The authors declare that they have no competing interests.

\bibliographystyle{elsarticle-num}
\bibliography{Bibliography}

\appendix

\end{document}